\documentclass[11pt]{article}  

\usepackage{hyperref}

\usepackage{amsfonts}
\usepackage{graphicx,subfig}
\usepackage{amsmath,amsthm}
\usepackage{wrapfig}

\usepackage{epsfig}
\usepackage{amsmath}
\usepackage{amssymb}
\usepackage{psfrag}
\usepackage[absolute]{textpos}

\newtheorem{thm}{Theorem}
\newtheorem{lem}{Lemma}

\def\x{{\bf x}}
\def\y{{\bf y}}
\def\w{{\bf w}}
\def\v{{\bf v}}

\newcommand{\el}{{\ell_1}}
\newcommand{\lt}{{\ell_2}}
\newcommand{\la}{\lambda}
\newcommand{\A}{\mathcal{A}}
\newcommand{\N}{\mathcal{N}}
\newcommand{\PS}{\mathbb{S}}
\newcommand{\R}{\mathbb{R}}
\newcommand{\G}{\mathcal{G}}
\newcommand{\HH}{\mathcal{D}}
\newcommand{\E}{\mathbb{E}}
\newcommand{\X} {{\bf{X}}}
\newcommand{\s}{\star}
\newcommand{\eps}{\epsilon}
\newcommand{\sig}{\sigma}

\newcommand{\h}{{\bf h}}
\newcommand{\z}{{\bf z}}
\newcommand{\dd}{{\bf d}}

\newcommand{\beq}{\begin{equation}}
\newcommand{\eeq}{\end{equation}}
\newcommand{\bea}{\begin{eqnarray}}
\newcommand{\eea}{\end{eqnarray}}

\newcommand{\Prob}{\ensuremath{\mathbb{P}}}

\long\def\symbolfootnote[#1]#2{\begingroup%
\def\thefootnote{\fnsymbol{footnote}}\footnote[#1]{#2}\endgroup}

\newtheorem*{WeakDef}{Weak recovery threshold}
\newtheorem*{SecDef}{Sectional recovery threshold}
\newtheorem*{StrDef}{Strong recovery threshold}
\newtheorem*{PSDWeak}{PSD Weak threshold}
\newtheorem*{PSDStr}{PSD Strong threshold}
\newtheorem*{UniWeak}{Uniqueness Weak threshold}
\newtheorem*{UniStr}{Uniqueness Strong threshold}

\begin{document}

\title{New Null Space Results and Recovery Thresholds for Matrix Rank Minimization}

\author{Samet Oymak,\, Babak Hassibi \\
  California Institute of Technology \\ Email: \{soymak,hassibi\}@caltech.edu
\thanks{This work was supported in part by the National Science Foundation under grants CCF-0729203, CNS-0932428 and CCF-1018927, by the Office of Naval Research under the MURI grant N00014-08-1-0747, and by Caltech's Lee Center for Advanced Networking.}
}
\maketitle
\begin{abstract}
Nuclear norm minimization (NNM) has recently gained significant attention for its use in rank minimization problems. Similar to compressed sensing, using null space characterizations, recovery thresholds for NNM have been studied in \cite{arxiv,Recht_Xu_Hassibi}. However simulations show that the thresholds are far from optimal, especially in the low rank region. In this paper we apply the recent analysis of Stojnic for compressed sensing \cite{mihailo} to the null space conditions of NNM. The resulting thresholds are significantly better and in particular our weak threshold appears to match with simulation results. Further our curves suggest for any rank growing linearly with matrix size $n$ we need only three times of oversampling (the model complexity) for weak recovery. Similar to \cite{arxiv} we analyze the conditions for weak, sectional and strong thresholds. Additionally a separate analysis is given for special case of positive semidefinite matrices. We conclude by discussing simulation results and future research directions.
\end{abstract}

\section{Introduction}

Rank minimization (RM) addresses the recovery of a low rank matrix from a set of linear measurements that project the matrix onto a lower dimensional space.  The problem has gained extensive attention in the past few years, due to the promising applicability in many practical problems \cite{Fazel}. Suppose that $X_0$ is a low rank matrix of size $n_1 \times n_2$ and let $\text{rank}(X)=r$. Further let $\A: \R^{n_1\times n_2} \rightarrow \R^{m}$ be a linear measurement operator. Given the measurements $y_0=\A(X_0)$, the problem is to recover $X_0$, with the knowledge of the fact that it is low rank. Provided that $X_0$ is the solution with lowest rank, this problem can be formulated with the following minimization program.
\begin{eqnarray}
&&\min~\text{rank}(X) \label{eq:rank min}\\
&&\mbox{subject to}\nonumber\\
 &&\hspace{20 pt}\A(X)=y_0, \nonumber
 \end{eqnarray}

\noindent The $\text{rank}(\cdot)$ function is non-convex, and it turns out that (\ref{eq:rank min}) is NP hard and cannot be solved efficiently. Fazel et al. suggested replacing the rank with the nuclear norm heuristic as the closest convex relaxation~\cite{Fazel}. The resulting convex optimization program is called nuclear norm minimization and is as follows.
\begin{eqnarray}
&&\min~\| X \|_\s \label{eq: nuc norm min}\\
&&\mbox{subject to}\nonumber\\
 &&\hspace{20 pt}\A(X)=\A(X_0), \nonumber
 \end{eqnarray}

\noindent where $\|\cdot\|_{\s}$ refers to the nuclear norm of its argument, i.e., the sum of the singular values. (\ref{eq: nuc norm min}) can be written as a semi-definite program (SDP) and thus be solved in polynomial time. Recent works have studied the sufficient conditions under which (\ref{eq: nuc norm min}) will recover $X_0$ (i.e. $X_0$ is unique minimizer of (\ref{eq: nuc norm min})). In \cite{Recht_Fazel_RIP} it is shown that, similar to compressed sensing, Restricted Isometry Property (RIP) is a sufficient condition for the success of (\ref{eq: nuc norm min}) and $O(rn_1n_2(n_1+n_2)log(n_1n_2))$ measurement is enough for guaranteeing RIP with high probability.  In \cite{candes_last}, Candes extended these results and showed that a minimal sampling of $O(rn)$ is in fact enough to have RIP and hence recovery. In later works \cite{Recht_Xu_Hassibi,arxiv}, necessary and sufficient null space conditions are derived and were analyzed for Gaussian measurement operators, i.e., operators where the entries are i.i.d. Gaussian, leading to thresholds for the success of (\ref{eq: nuc norm min}). These thresholds establish explicit relationships between the problem parameters, as opposed to the order-wise relationships that result from RIP techniques. However these results are far from being optimal in the low rank regime which necessitates a new approach to be taken. In particular, if the matrix size is $n\times n$ and the rank of the matrix to be recovered is $\beta n$ then even if $\beta>0$ is very small, they require a minimum sampling of $(1-\frac{64}{9\pi^2})n^2$ for success. In this paper, we come up with a novel null space analysis for the rank minimization problem and we find significantly better thresholds than the results of \cite{Recht_Xu_Hassibi, arxiv}. Although the analysis is novel for the rank minimization problem, we basically follow the analysis developed for compressed sensing by Stojnic in \cite{mihailo} which is based on a seminal result of Gordon \cite{Gordan2}. In addition to the analysis of general matrices, we give a separate analysis for positive semidefinite matrices which resemble nonnegative vectors in compressed sensing. We also consider the case of unique positive semidefinite solutions, which was recently analyzed by Xu in \cite{weiyu}.

We extensively use the results of \cite{mihailo}. Basically, we slightly modify Lemmas 2, 5, 7 of \cite{mihailo} and use null space conditions for the NNM problem. The strength of this analysis comes from the facts that the analysis is more accessible and that the weak threshold of \cite{mihailo} matches the exact threshold of \cite{Donoho2}. In fact, while it is not at all clear how to extend the analysis of \cite{Donoho2} from compressed sensing to NNM, it is relatively straightforward to do so for \cite{mihailo}. Our simulation results also indicate that our thresholds for the NNM problem are seemingly tight. This is perhaps not surprising since, as we shall see, the null space conditions for NNM and compressed sensing are very similar.

\section{Basic Definitions and Notations}

Denote identity matrix of size $n\times n$ by $I_n$. We call $U\in\R^{m\times n}$ {\bf{partial unitary}} if columns of $U$ form an orthonormal set i.e. $U^TU=I_n$. Clearly we need $n\leq m$ for $U$ to be partial unitary. Also for a partial unitary $U$, let $\bar{U}$ denote an arbitrary partial unitary of size $m\times (m-n)$ so that $[U~\bar{U}]$ is a unitary matrix (i.e. columns are complete orthonormal basis of $\R^m$). 

For a matrix $X\in \R^{m\times n}$, we denote the singular values by $\sig_1(X)\geq \sig_2(X)\geq \dots \geq \sig_q(X)$ where $q=\min(m,n)$. The (skinny) singular value decomposition (SVD) of $X$ is shown as $X=U_X \Sigma_X V_X^T$ where $U_X\in\R^{m\times r}$, $\Sigma_X\in\R^{r\times r}$ and $V_X\in\R^{n\times r}$, where $r=\text{rank}(X)$. Note that $U_X$,$V_X$ are partial unitary and $\Sigma_X$ is positive, diagonal and full rank. Also let $\Sigma(X)$ denote vector of increasingly ordered singular values of $X$ i.e. $\Sigma(X)=[\sig_q(X)~\dots~\sig_1(X)]^T$.

The $\text{Ky-Fan}~k$ norm of $X$ denoted by $\| X\|_k$  is defined as $\| X\|_k=\sum_{i=1}^k \sig_i(X)$.
When $k=\min(m,n)$ it is called the nuclear norm, i.e. $\| X\|_{\s}$, and when $k=1$ it is equivalent to the spectral norm denoted by $\| X\|$. Also Frobenius norm is denoted by $\| X\|_F=\sqrt{\left<X,X\right>}=\sqrt{\sum_{i=1}^n\sig_i^2(X)}$. Note that we always have:
\beq
\|X\|_k=\sum_{i=1}^k\sig_i(X)\leq \sqrt{\sum_{i=1}^k 1\sum_{i=1}^k\sig_i^2(X)}\leq \sqrt{k}\|X\|_F
\eeq

For a linear operator $\A(\cdot)$ acting on a linear space, we denote the null space of $\A$ by $\N(\A)$, i.e. $W\in\N(\A)$ iff $\A(W)=0$. We denote by $\G(d_1,d_2)$ the ensemble of real $d_1\times d_2$ matrices in which the entries are i.i.d. $\N(0,1)$ (zero-mean, unit variance Gaussian).

It is a well known fact that normalized singular values of a square matrix with i.i.d. Gaussian entries have quarter circle distribution asymptotically \cite{marpas}. In other words the histogram of singular values (normalized by $1/\sqrt{n}$) converges to the function
\begin{align}
\phi(x)=\frac{\sqrt{4-x^2}}{\pi}~~~~0\leq x\leq 2
\end{align}
Similarly, the distribution of the squares of the singular values (normalized by $1/n$) converges to the well known Marcenko-Pastur distribution \cite{marpas}. Note that this is nothing but the distribution of the eigenvalues of $X^TX$ where $X$ is a square matrix drawn from $\G(n,n)$,
\begin{align}
\phi_2(x)=\frac{\sqrt{4x-x^2}}{2\pi x}~~~~0\leq x\leq 4
\end{align}
Let $F(x)$ be the cumulative distribution function of $\phi(x)$ i.e.,
\begin{align}
\label{eq:cdfs}
F(x)=\int_{0}^{x} \phi(t) dt
\end{align}
Let $0\leq\beta\leq 1$. We define $\gamma(\beta)$ to be the asymptotic normalized expected value of the $\text{Ky-Fan}~\beta n$ norm of a matrix drawn from $\G(n,n)$, i.e.:
\begin{equation}
\gamma(\beta):=\lim_{n\rightarrow\infty}\frac{\E[\|X\|_{\beta n}]}{n^{3/2}}=\int_{F^{-1}(1-\beta)}^{2} x\phi(x) dx
\end{equation}
Similarly define $\gamma_2(\beta)$ to be the asymptotic normalized expected value of the $\text{Ky-Fan}~\beta n$ norm of a matrix $X^TX$ where $X$ is drawn from $\G(n,n)$:
\begin{align}
\gamma_2(\beta)=\lim_{n\rightarrow\infty}\E\left[\frac{[\|X^TX\|_{\beta n}]}{n^2}\right]=\lim_{n\rightarrow\infty} \E\left[\frac{\sum_{i=1}^{\beta n}\sig_i(X)^2}{n^2}\right]
=\int_{F^{-1}(1-\beta)}^2 x^2\phi(x)dx
\end{align}

\noindent Note that these limits exist and $\gamma(\beta),\gamma_2(\beta)$ is well defined \cite{Random Matrix}.

A function $f:\R^n\rightarrow \R$ is called $L$-Lipschitz if for all $x,y$ we have: $|f(x)-f(y)|\leq L\|x-y\|_{\ell_2}$

We say an orthogonal projection pair $\{P,Q\}$ is a support of the matrix $X$ if $X=PXQ^T$. In particular $\{P_X,Q_X\}$ is the unique support of the matrix $X$, if $P_X$ and $Q_X$ are orthogonal projectors with $\text{rank}(P_X)=\text{rank}(Q_X)=\text{rank}(X)$ such that $X=P_XXQ_X^T$. In other words, $P_X = U_XU_X^T$ and $Q_X = V_XV_X^T$.

We say $\A:\R^{n_1\times n_2}\rightarrow\R^m$ is a random Gaussian measurement operator if the $i^{\text{th}}$ measurement is $y_i=\A(X)_i=\text{trace}(G_i^TX)$ where $\{G_i\}_{i=1}^m$'s are i.i.d. matrices drawn from $\G(n_1,n_2)$ for all $1\leq i\leq m$. Note that this is equivalent to
$y_i=\text{vec}(G_i)^T\text{vec}(X)$ where $\text{vec}(X)$ is obtained by putting columns of $X$ on top of each other to get a vector of size $n_1 n_2\times 1$.

Model complexity is defined as the number of degrees of freedom of the matrix. For a matrix of size $\alpha n\times n$ and rank $\alpha\beta n$ model complexity is $\alpha\beta(1+\alpha-\alpha\beta+o(1))n^2$. Then we define normalized model complexity to be $\theta=\beta(1+\alpha-\alpha\beta)$.

Finally let $\succeq$ denote "greater than" in partially ordered sets. In particular if $A,B$ are Hermitian matrices then $A\succeq B\iff A-B~\text{is positive semidefinite}$. Similarly for a given two vectors $u,v$ we write $u\succeq v\iff u_i\geq v_i~\forall~i$.

\section{Key Lemmas to be Used}
In this section, we state several lemmas that we will make use of later. Proofs that are omitted can be found in the given references.

For Lemmas (\ref{lem:neumann}), (\ref{lem:singval}), (\ref{lem:simple}), let $X,Y,Z\in \R^{m\times n}$ with $m\leq n$.
\begin{lem}
\label{lem:neumann}
\begin{eqnarray}
tr(X^TY)\leq \sum_{i=1}^m \sig_i(X)\sig_i(Y)=\Sigma(X)^T\Sigma(Y)
\end{eqnarray}
\begin{proof}
Can be found in \cite{mirsky}.
\end{proof}
In case of vectors (i.e. matrices are diagonal) we have the following simple extension: Let $\x,\y\in\R^m$ be vectors. Let $x_{[i]}$ be $i$'th largest value of vector $|\x|$ (i.e. $|\x|_i=|\x_i|$) then
\begin{equation}
\left<\x,\y\right>\leq \sum_{i=1}^m x_{[i]} y_{[i]}
\end{equation}
\end{lem}

\begin{lem}
\label{lem:singval}
Let $Z=X-Y$. Let $s_i(X,Y)=|\sig_i(X)-\sig_i(Y)|$ and let $s_{[1]}(X,Y)\geq s_{[2]}(X,Y)\geq \dots \geq s_{[m]}(X,Y)$ be a decreasingly ordered arrangement of $\{s_i(X,Y)\}_{i=1}^m$. Then we have the following inequality:
\begin{equation}
\forall~m\geq k\geq 1:~~\sum_{i=1}^{k}s_{[i]}(X,Y) \leq \sum_{i=1}^{k} \sig_i(Z)=\| Z\|_k
\end{equation}
In particular we have:
\begin{eqnarray}
\sum_{i=1}^{m}|\sig_i(X)-\sig_i(Y)|\leq \| Z \|_{\s}
\end{eqnarray}
\begin{proof}
Proof can be found in \cite{Horn,Bhatia}
\end{proof}
\end{lem}

\begin{lem}
\label{lem:simple}
If matrix $X=\begin{bmatrix} X_{11} & X_{12} \\ X_{21} & X_{22}\end{bmatrix}$ then we have:
\begin{equation}
\| X\|_\s\geq \| X_{11}\|_\s+\| X_{22}\|_\s
\end{equation}
\begin{proof}
Proof can be found in \cite{Recht_Xu_Hassibi}.
\end{proof}
Similarly, we have the following obvious inequality when $X$ is square ($m=n$):
\beq
\| X\|_\s\geq \text{trace}(X)
\eeq
\begin{proof}
Dual norm of the nuclear norm is the spectral norm \cite{Horn}. Remember that $I_m$ is identity. Then:
\beq
\|X\|_\s=\sup_{\|Y\|=1} \left<X,Y\right>\geq \left<X,I_m\right>=\text{trace}(X)
\eeq
\end{proof}
\end{lem}

\begin{thm} \label{thm: ETM} (Escape through a mesh, \cite{Gordan2}) Let $S$ be a subset of the unit Euclidean sphere $S^{n-1}$ in $\R^n$. Let $Y$ be a random $(n - m)$-dimensional subspace of $\R^n$, distributed uniformly in the Grassmanian with respect to Haar measure. Let
\begin{equation}
\label{eq: ETMeq}
\omega(S)=\E \sup_{\w\in S} (\h^T \w)
\end{equation}
where $\h$ is a column vector drawn from $\G(n,1)$. Then if $\omega(S)<\sqrt{m}-\frac{1}{4\sqrt{m}}$ we have:
\begin{equation}
\Prob(Y \cap S=\emptyset)>1-3.5\exp(-\frac{\left(\sqrt{m}-\frac{1}{4\sqrt{m}}-\omega(S)\right)^2}{18}
\end{equation}
\end{thm}

\begin{lem}
\label{lem:lipsc1}
For all $1\leq k\leq n$, $\sig_k(X)$ is a $1\text{-Lipschitz}$ function of $X$.
\begin{proof}
Let $X,\hat{X}$ be such that $\tilde{X}=\hat{X}-X$ and $\|\tilde{X}\|_F\leq 1$. But then from Lemma (\ref{lem:singval}) we have:
\begin{equation}
1\geq \|\tilde{X}\|_F\geq \sig_1(\tilde{X})\geq s_{[1]}(X,\hat{X})\geq |\sig_k(X)-\sig_k(\hat{X})|
\end{equation}
\end{proof}
\end{lem}

\begin{lem}
\label{lem:gauslip}
(from \cite{Recht_Xu_Hassibi,talagrand}) Let $x$ be drawn from $\G(n,1)$ and $f : \R^n \rightarrow \R$ be a function with Lipschitz constant $L$ then we have the following concentration inequality
\begin{equation}
\Prob(|f(x)-\E f(x)|\geq t)\leq 2\exp(-\frac{t^2}{2L^2})
\end{equation}
\end{lem}

For analyzing positive semidefinite matrices, we will introduce some more definitions and lemmas later on.

\section{Thresholds for Square Matrices}

In the following section, we'll give and analyze strong, sectional and weak null space conditions for square matrices ($\R^{n\times n}$). With minor modifications, one can obtain the equivalent results for rectangular matrices ($\R^{\alpha n\times n}$).


\subsection{Strong Threshold}
\begin{StrDef}
Let $\A:\R^{n\times n}\rightarrow \R^m$ be a random Gaussian operator. We define $\beta$ ($0\leq\beta\leq 1$) to be the strong recovery threshold if with high probability $\A$ satisfies the following property:

Any matrix $X$ with rank at most $\beta n$ can be recovered from measurements $\A(X)$ via (\ref{eq: nuc norm min}).
\end{StrDef}
\begin{lem}
\label{lem:strong_null}
Using (\ref{eq: nuc norm min}) one can recover all matrices $X$ of rank at most $r$ if and only if for all $W\in \N(\A)$ we have
\begin{equation}
\label{eq: str null}
2\| W\|_r<\| W\|_\s
\end{equation}
\begin{proof}
If (\ref{eq: str null}) holds then using Lemma (\ref{lem:singval}) and the fact that $\sig_i(X)=0~~\forall~i>r$ for any $W$ we have
\begin{eqnarray}
\| X+W\|_\s\geq \sum_{i=1}^n |\sig_i(X)-\sig_i(W)|&\geq& \sum_{i=1}^r (\sig_i(X)-\sig_i(W)) +\sum_{i=r+1}^n \sig_i(W)\\
&\geq&\| X\|_\s +\| W\|_\s -2\|W\|_r > \| X\|_\s
\end{eqnarray}
Hence $X$ is unique minimizer of (\ref{eq: nuc norm min}).
Conversely if (\ref{eq: str null}) doesn't hold for some $W$ then choose $X=-W_k$ where $W_k$ is the matrix induced by setting all but largest $r$ singular values of $W$ to 0. Then we get: $\| X+W\| =\sum_{i=r+1}^n \sig_i(W)\geq  \sum_{i=1}^r \sig_i(W)=\sum_{i=1}^r \sig_i(X)=\|X\|_\s$. Finally we find $\text{rank}(X)\leq r$ but $X$ is not the unique minimizer.
\end{proof}
\end{lem}

Now we can start analyzing the strong null space condition for the NNM problem. $\A$ is a random Gaussian operator and we'll analyze the linear regime where $m=\mu n^2$ and $r=\beta n$. Our aim is to determine the least $\mu$ ($1\geq\mu\geq0$) so that $\beta$ is a strong threshold for $\A$. Similar to compressed sensing the null space of $\A$ is an $n^2-m$ dimensional random subspace of $\R^{n^2}$ distributed uniformly in the Grassmanian w.r.t. Haar measure. This can also be viewed as the span of $M=(1-\mu) n^2$ matrices $\{ G_i\}_{i=1}^{M}$ drawn i.i.d. from $\G(n,n)$. Then similar to \cite{mihailo} we have established the necessary framework.

Let $S_s$ be the set of all matrices such that $\| W\|_\s \leq 2\| W\|_r$ and $\| W\|_F=1$. We need to make sure the null space of $\A$ has no intersection with $S_s$. We will first upper bound (\ref{eq: ETMeq}) in Theorem \ref{thm: ETM} then choose $m$ (and $\mu$) respectively.

As a first step, given a fixed $H\in \R^{n\times n}$ we'll calculate an upper bound on $f(H,S_s)=\sup_{W\in S_s} vec(H)^T vec(W)=\sup_{W_\in S_s} \left< H,W\right> $. Note that from Lemma \ref{lem:neumann} we have:
\begin{equation}
\label{eq:rasgele1}
f(H,S_s)=\sup_{W \in S_s} \left< H,W\right>\leq \sup_{W \in S_s} \Sigma(H)^T\Sigma(W)
\end{equation}
The careful reader will notice that actually we have equality in (\ref{eq:rasgele1}) because the set $S_s$ is unitarily invariant hence any value we can get on the right hand side, we can also get on the left hand side by aligning the singular vectors of $H$ and $W$. Let $\h=\Sigma(H)$, $\w=\Sigma(W)$. Note that $\h,\w\succeq 0$. Then since $\|\w\|_{\ell_2}=\|W\|_F=1$ and $\sum_{i=1}^{n-r}w_i\leq \sum_{i=n-r+1}^nw_i$ any $W\in S_s$, we need to solve the following optimization problem given $H$:
\begin{eqnarray}
\label{eq:optim1}
&&\max_\y~ \h^T \y\\
&&\nonumber\text{subject to}\\
&&\nonumber~~~~\y\succeq 0\\
&&\nonumber~~~~\sum_{i=n-r+1}^n y_i \geq \sum_{i=1}^{n-r} y_i\\
&&\nonumber~~~~\| \y\|_{\ell_2}\leq1
\end{eqnarray}
Clearly the right hand side of (\ref{eq:rasgele1}) and the result of (\ref{eq:optim1}) is same because $\h^T \y$ will be maximized when $\{y_i\}_{i=1}^{n}$ are sorted increasingly due to Lemma \ref{lem:neumann}.

Note that (\ref{eq:optim1}) is exactly the same as (10) of \cite{mihailo}. Then we can use (22), (29) of \cite{mihailo} directly to get:
\begin{lem}
\label{lem:strpes0}
If $\h^T\z >0$ then
\begin{equation}
\label{eq:str1}
f(H,S_s)\leq \sqrt{\sum_{i=c+1}^n h_i^2 -\frac{((\h^T \z)-\sum_{i=1}^c h_i)^2}{n-c}}
\end{equation}
where $\z \in \R^n$ such that $z_i=1~~\forall~~1\leq i\leq n-r$ and $z_i=-1~~\forall~~n-r+1\leq i\leq n$ and $0\leq c\leq n-r$ such that $(\h^T \z)-\sum_{i=1}^c h_i\geq (n-c)h_c$. As long as $\h^T\z> 0$ we can find such $c\geq 0$. In addition, in order to minimize right hand side of (\ref{eq:str1}), one should choose largest such $c$.

In case of $\h^T\z\leq 0$, the following is the obvious upper bound from Cauchy-Schwarz and the fact that $\| W\|_F=1$
\begin{equation}
f(H,S_s)\leq \| \h\|_{\ell_2}=\sqrt{\sum_{i=1}^n h_i^2}
\end{equation}
\end{lem}

Similar to \cite{mihailo}, for the escape through a mesh (ETM) analysis, using Lemma \ref{lem:strpes0}, we'll consider the following worse upper bound:
\begin{lem}
\label{lem:strpes1}
Let $\z$ be defined same as in Lemma \ref{lem:strpes0}. Let $H$ be chosen from $\G(n,n)$ and let $\h=\Sigma(H)$ and $f(H,S_s)=\sup_{W\in S_s} \left<H,W\right>$. Then we have: $f(H,S_s)\leq B_s$ where
\begin{align}
&\nonumber B_s= \| \h\|_{\ell_2}~~~\text{if }g(H,c_s)\leq 0\\
&\nonumber B_s=\sqrt{\sum_{i=c_s+1}^n h_i^2 -\frac{((\h^T \z)-\sum_{i=1}^{c_s} h_i)^2}{n-c_s}}~~~\text{else}
\end{align}
where $g(H,c)=\frac{(\h^T \z)-\sum_{i=1}^c h_i}{n-c}-h_c$ and $c_s=\delta_s n$ is a $c\leq n-r$ such that
\begin{align}
\label{eq:csst}
\nonumber c_s&=0~~~\text{ if}~\E[\h^T\z]\leq 0\\
c_s&\text{ is solution of }~~~(1-\eps)\frac{\E [(\h^T \z)-\sum_{i=1}^c h_i]}{\sqrt{n}(n-c)}=F^{-1}\left(\frac{(1+\eps)c}{n}\right)\text{ else if }\E[\h^T\z]>0
\end{align}
where $\eps>0$ can be arbitrarily small. Note that $c_s$ is deterministic. Secondly one can observe that $c_s>0\iff \E[\h^T\z]>0\iff \E[\h^T\z-\sum_{i=1}^{c_s} h_i]>0$.
\end{lem}
Here $F(\cdot)$ is the c.d.f. of the quarter circle distribution previously defined in (\ref{eq:cdfs}).

\subsubsection{Probabilistic Analysis of $\E[B_s]$}
The matrix $H$ is drawn from $\G(n,n)$ and $\E[B_s]\geq \E[f(H,S_s)]$. In the following discussion, we'll focus on the case $\E[\h^T\z]>0$ and we'll declare failure (no recovery) else. This is reasonable since our approach will eventually lead to $\mu=1$ in case of $\E[\h^T\z]\leq 0$. The reason is that, with high probability we'll have $g(H,c_s)\leq 0$ and this will result in $\E[B_s]\approx\E[\|H\|_F]$ which is the worst upper bound. 

Then, we'll basically argue that whenever $\E[\h^T\z]>0$, asymptotically with probability one, we'll have $g(H,c_s)>0$. Next, we'll show that contribution of the region $g(H,c_s)\leq0$ to the expectation of $B_s$ asymptotically converges to $0$.

From the union bound, we have:
\begin{eqnarray}
\Prob(g(H,c_s)\leq 0)\leq \Prob(\h^T \z-\sum_{i=1}^{c_s}h_i \leq (1-\eps)\E[(\h^T \z)-\sum_{i=1}^{c_s} h_i])+\Prob(h_{c_s}\geq \sqrt{n}F^{-1}\left(\frac{(1+\eps)c_s}{n}\right))
\end{eqnarray}

We'll analyze the two components separately. Note that $\h^T\z-\sum_{i=1}^{c_s}h_i$ is a function of singular values which is actually a Lipschitz function of the random matrix $H$ as we'll argue in the following lemma.
\begin{lem}
Let $H\in\R^{n\times n}$ and let $\h=\Sigma(H)$ and $\z$ is as defined previously. Then:
\begin{equation}
f(H)=\h^T\z-\sum_{i=1}^{c_s}h_i
\end{equation}
is $\sqrt{n-c_s}$ Lipschitz function of $H$.
\begin{proof}
Let $H,\hat{H},\tilde{H}$ be such that $\tilde{H}=H-\hat{H}$. From Lemma (\ref{lem:singval}) we have:
\begin{align}
\| \tilde{H}\|_{n-c_s}\geq \sum_{i=1}^{n-c_s} |\sig_i(H)-\sig_i(\hat{H})|&\geq  |\sum_{i=1}^r (\sig_i(H)-\sig_i(\hat{H}))|+|\sum_{i=r+1}^{n-c_s} (\sig_i(\hat{H})-\sig_i(H))|\\
&\geq |\h^T\z -\hat{\h}^T\z| = |f(H)-f(\hat{H})|
\end{align}
On the other hand we have: $\| \tilde{H}\|_{n-c_s} \leq \sqrt{n-c_s}\|\tilde{H}\|_F$ which implies $|f(H)-f(\hat{H})|\leq \sqrt{n-c_s}\| \tilde{H}\|_F$ finishing the proof.
\end{proof}
\end{lem}
Now, using the fact that $H$ is i.i.d. Gaussian and $\h$ is the vector of singular values of $H$, we have $\E(\h^T\z-\sum_{i=1}^{c_s}h_i)=(\gamma(1)-2\gamma(\beta)-(\gamma(1)-\gamma(1-{\delta_s}))+o(1))n^{3/2}$ hence from Lemma \ref{lem:gauslip} and from the fact that $H$ is i.i.d. Gaussian, we have:
\begin{equation}
P_1:=\Prob(\h^T\z-\sum_{i=1}^{c_s}h_i \leq (1-\eps)\E[(\h^T \z)-\sum_{i=1}^{c_s} h_i])\leq \exp(-\frac{\eps^2(\gamma(1-\delta_s)-2\gamma(\beta))+o(1))^2n^2}{2(1-\delta_s)})
\end{equation}
if $\E[\h^T\z]> 0$ (which is equivalent to $\E[\h^T\z-\sum_{i=1}^{c_s}h_i]> 0$ and $\delta_s>0$).

Similarly from the quarter circle law we have $\E(h_c)=(F^{-1}(c/n)+o(1))\sqrt{n}$. Using Lemmas \ref{lem:gauslip}, \ref{lem:lipsc1} we can find:
\begin{align}
P_2:=\Prob(h_{c_s}\geq \sqrt{n}F^{-1}\left(\frac{(1+\eps)c_s}{n}\right))\leq \exp(-\frac{n}{2} \left(F^{-1}\left(\frac{(1+\eps)c_s}{n}\right)-F^{-1}\left(\frac{c_s}{n}\right)+o(1)\right)^2)
\end{align}
In particular we always have $F^{-1}((1+\eps)\beta)-F^{-1}(\beta)\geq \frac{\pi\eps\beta}{2}$ for any $\eps>0$, $1>\beta\geq0$ (because $F(x)\leq \frac{2}{\pi}$ for $0\leq x\leq 2$). Hence $P_2$ converges to $0$ exponentially fast. One can actually show $P_2\leq \exp(-O(n^2))$ instead of $\exp(-O(n))$ however this won't affect the results.

Then since $\delta_s>0$: $P(g(H,c_s)\leq 0)\leq P_1+P_2\leq \exp(-\frac{n}{8}(\pi\eps\delta_s+o(1))^2)$. It remains to upper bound $\E(B_s)$ as follows:
\begin{equation}
\label{eq:str2}
\E[B_s]\leq \int_{g(H,c_s)\leq 0} \| \h\|_{\ell_2}p(H)dH +\int_H \sqrt{\sum_{i=c_s+1}^n h_i^2 -\frac{((\h^T \z)-\sum_{i=1}^{c_s} h_i)^2}{n-c_s}}p(H)dH
\end{equation}
Note that $g(H,c)$ is linear function of $h$ (hence $H$) so if $g(H,c)\leq 0\iff g(aH,c)\leq 0$ for any $a\geq 0$. In other words similar to the discussion in \cite{mihailo} for any value of $a=\| H\|_{F}$, the fraction of the region $g(H,c)\leq 0$ on the sphere of radius $a$ will be constant. On the other hand since $H$ is iid Gaussian, the probability distribution of $H$ is just a function of $\|H\|_F$ i.e. $p(H=H_0)=f(\| H_0\|_F)=(2\pi)^{-n^2/2}\exp(-\frac{1}{2}\|H_0\|_F^2)$ for any matrix $H_0\in\R^{n\times n}$. As a result:
\begin{equation}
\int_{g(H,c_s)\leq 0,\|H\|_F=a}dH=C_0\int_{|H\|_F=a}dH=C_0S_a
\end{equation}
where $S_a$ is the area of a sphere in $\R^{n\times n}$ with radius $a$. Hence
\begin{eqnarray}
P(g(H,c_s)\leq 0)&=&\int_{a\geq 0}\int_{g(H,c_s)\leq 0,\|H\|_F=a}p(H)dHda=\int_{a\geq 0}\int_{g(H,c_s)\leq 0,\|H\|_F=a}f(a)dHda\\
&=&C_0\int_{a\geq 0} f(a)S_a da=C_0
\end{eqnarray}

Using the exact same argument:
\begin{eqnarray}
\label{eq:goeszero}
\nonumber\int_{g(H,c_s)\leq 0} \| H\|_{F}p(H)dH &=&\int_{a=0}^{\infty}\int_{g(H,c_s)\leq 0,\|H\|_F=a} \| H\|_Fp(H)dHda\\
\nonumber&=&\int_{a=0}^\infty \int_{g(H,c_s)\leq 0,\|H\|_F=a} af(a)dHda\\
\nonumber&=&\int_{a=0}^\infty af(a)C_0 S_a=P(g(H,c_s)\leq 0)\E(\| H\|_F)\\
&\leq&\exp(-\frac{n}{8}(\pi\eps\delta_s+o(1))^2)n
\end{eqnarray}
The last term clearly goes to zero for large $n$. Then we need to calculate the second part which is:
\begin{eqnarray}
\int_H \sqrt{\sum_{i=c_s+1}^n h_i^2 -\frac{((\h^T \z)-\sum_{i=1}^{c_s} h_i)^2}{n-c_s}}p(H)dH&=&\E( \sqrt{\sum_{i=c_s+1}^n h_i^2 -\frac{((\h^T \z)-\sum_{i=1}^{c_s} h_i)^2}{n-c_s}})\\
&\leq&\sqrt{\E(\sum_{i=c_s+1}^n h_i^2 -\frac{((\h^T \z)-\sum_{i=1}^{c_s} h_i)^2}{n-c_s})}
\end{eqnarray}
The last inequality is due to the following Cauchy-Schwarz. For a random variable (R.V.) $\X\geq0$
\begin{equation}
\E(\X)=\int_x x p(x)dx\int_x p(x)dx\geq \left(\int_x \sqrt{xp(x)^2}dx\right)^2=\E(\sqrt{\X})^2
\end{equation}
Note that for large $n$ and fixed $c_s=\delta_s n$ and $r=\beta n$ we have
\begin{equation}
\label{eq:eqstr}
\E(\sum_{i=c_s+1}^n h_i^2 -\frac{((\h^T \z)-\sum_{i=1}^{c_s} h_i)^2}{n-c_s})=\left(\gamma_2(1-\delta_s)-\frac{(\gamma(1-\delta_s)-2\gamma(\beta))^2}{1-\delta_s}+o(1)\right)n^2
\end{equation}
Then combining (\ref{eq:str2}) and (\ref{eq:goeszero}), it follows that (\ref{eq:eqstr}) gives an upper bound for $\E[B_s]^2$ and thereby $\E[f(H,S_s)]^2$. To be able to calculate the required number of measurements we need to find $\delta_s$ and substitute in (\ref{eq:eqstr}) because (\ref{eq:eqstr}) will also be an upper bound on the minimum $m$ asymptotically.

If we consider (\ref{eq:csst}), asymptotically $\delta_s$ will be solution of:
\begin{equation}
\label{eq:strdelta}
(1-\eps)\frac{\gamma(1-\delta_s)-2\gamma(\beta)}{1-\delta_s}=F^{-1}((1+\eps)\delta_s)
\end{equation}
Then we can substitute this $\delta_s$ in (\ref{eq:eqstr}) to solve for $m$ (and $\mu$). Using Theorem \ref{thm: ETM} and (\ref{eq:strdelta}) we find:
\begin{thm}
\label{thm:strfinal}
If $\gamma(1)-2\gamma(\beta)\leq0$ then $\mu=1$. Otherwise:
\begin{equation}
\label{eq:str3}
\mu>\gamma_2(1-\delta_s)-\frac{(\gamma(1-\delta_s)-2\gamma(\beta))^2}{1-\delta_s}
\end{equation}
is sufficient sampling rate for $\beta$ to be strong threshold of random Gaussian operator $\A:\R^{n\times n}\rightarrow\R^{\mu n^2}$. Here $\delta_s$ is solution of:
\begin{equation}
(1-\eps)\frac{\gamma(1-\delta_s)-2\gamma(\beta)}{1-\delta_s}=F^{-1}((1+\eps)\delta_s)
\end{equation}
\end{thm}
In order to get the smallest $\mu$ we let $\eps\rightarrow 0$. Numerical calculations give the strong threshold in Figure \ref{fig:real_thresholds}. Obviously we found and plotted the least $\mu$ for a given $\beta$ (i.e. equality in (\ref{eq:str3})).\\ 

Next we define and analyze sectional threshold.

\subsection{Sectional Threshold}

\begin{SecDef}
Let $\A:\R^{n\times n}\rightarrow \R^m$ be a random Gaussian operator and let $\{P,Q\}$ be an arbitrary orthogonal projection pair with $\text{rank}(P)=\text{rank}(Q)= \beta n$. Then we say that $\beta$ ($0\leq\beta\leq 1$) is a sectional recovery threshold if with high probability $\A$ satisfies the following property:

Any matrix $X$ with support $\{P,Q\}$ can be recovered from measurements $\A(X)$ via (\ref{eq: nuc norm min}).
\end{SecDef}
Given a fixed $\beta$, our aim is to calculate the least $\mu$ such that $\beta$ is sectional threshold for a random Gaussian operator $\A:\R^{n\times n}\rightarrow \R^{\mu n^2}$.

\begin{lem}
\label{lem:sectional_null}
Given support $\{P,Q\}$ with $\text{rank}(P)=\text{rank}(Q)=r$ one can recover all matrices $X$ with this support using (\ref{eq: nuc norm min}) iff for all $W\in \N(\A)$ we have
\begin{equation}
\label{eq: sec null}
\|(I-P)W(I-Q^T)\|_\s>\| PWQ^T\|_\s
\end{equation}
\begin{proof}
Note that in a suitable basis induced by $\{P,Q\}$ we can write:
\begin{equation}
X=\begin{bmatrix}X_{11} & 0\\0& 0 \end{bmatrix},~~~~W=\begin{bmatrix}W_{11} & W_{12}\\W_{21}& W_{22} \end{bmatrix}
\end{equation}
where $\begin{bmatrix}W_{11} & 0\\0&0\end{bmatrix}=PWQ^T,~\begin{bmatrix}0 & W_{12}\\0&0\end{bmatrix}=PW(I-Q^T),~\begin{bmatrix}0 & 0\\W_{21}&0\end{bmatrix}=(I-P)WQ^T,~\begin{bmatrix}0 & 0\\0&W_{22}\end{bmatrix}=(I-P)W(I-Q^T)$. Now If (\ref{eq: sec null}) holds then using Lemma \ref{lem:simple} we immediately have for all $W\in\N(\A)$:
\begin{equation}
\|X+W\|_\s=\begin{bmatrix}X_{11}+W_{11} & W_{12}\\W_{21}& W_{22} \end{bmatrix}\geq \|X_{11}+W_{11}\|_\s+\|W_{22}\|_\s\geq \|X_{11}\|_\s-\|W_{11}\|_\s+\|W_{22}\|_\s>\|X_{11}\|_\s
\end{equation}
Hence $X$ is unique minimizer of (\ref{eq: nuc norm min}). In \cite{arxiv}, it was proven that (\ref{eq: sec null}) is tight in the sense that if there exists $W\in\N(\A)$ such that $\|(I-P)W(I-Q^T)\|_\s<\| PWQ^T\|_\s$ then we can find an $X$ with support $\{P,Q\}$ where $X$ is {\bf{not}} minimizer of  (\ref{eq: nuc norm min}).
\end{proof}
\end{lem}

Now we can start analyzing the sectional null space condition for the NNM problem. $\A$ is a random Gaussian operator and we'll analyze the linear regime where $m=\mu n^2$ and $r=\beta n$. Similar to compressed sensing, the null space of $\A$ is an $n^2-m$ dimensional random subspace of $\R^{n^2}$ distributed uniformly in the Grassmanian w.r.t. Haar measure. Then similar to \cite{mihailo} we have established the necessary framework.

Let $S_{sec}$ be the set of all matrices such that $\|(I-P)W(I-Q^T)\|_\s \leq \| PWQ^T\|_\s$ and $\| W\|_F=1$. We need to make sure, the null space has no intersection with $S_{sec}$. We will first upper bound (\ref{eq: ETMeq}) in Theorem \ref{thm: ETM}, then choose $m$ (and $\mu$) respectively. As discussed in \cite{arxiv}, without loss of generality we can assume $X=\begin{bmatrix} X_{11} & 0\\0&0\end{bmatrix}$ because $X$ can be transformed to this form with a unitary transformation (which depends only on $\{P,Q\}$) and since the null space is uniformly chosen (i.e. its basis is $n^2-m$ random matrices chosen iid from $\G(n,n)$) after this unitary transformation its distribution will still be uniform. The reason is that if $X$ is i.i.d. Gaussian matrix and $A,B$ are fixed unitary matrices then $AXB$ is still i.i.d. Gaussian. This further shows that the probability of successful recovery does not depend on $\{P,Q\}$ as long as $\beta$ is fixed. With this assumption $S_{sec}$ is the set of all matrices with $\|W_{22}\|_\s\leq \|W_{11}\|_\s$ and $\|W\|_F=1$. Observe that $W_{11}\in\R^{r\times r}$ and $W_{22}\in\R^{(n-r)\times (n-r)}$.

In the following we assume $2\times 2$ block matrices. Let $X_{ij}$ be $i$'th row and $j$'th column block of $X$. As a first step, given a fixed $H\in \R^{n\times n}$ we'll calculate an upper bound on $f(H,S_{sec})=\sup_{W\in S_{sec}} vec(H)^T vec(W)=\sup_{W_\in S_{sec}} \left< H,W\right> $. Note that: $ \left< H,W\right>=\sum_{i,j}  \left< H_{ij},W_{ij}\right>$ 

Further let $\h_1=\Sigma(H_{11}),~\h_2=\Sigma(H_{22}),~\w_1=\Sigma(W_{11}),~\w_2=\Sigma(W_{22})$. Also let $\h_3$ be increasingly sorted absolute values of entries of submatrices $H_{12},H_{21}$ and $\w_3$ is defined similarly. Finally let $x_{i,j}$ denote $j$'th entry of vector $\x_i$

From Lemma \ref{lem:neumann} we have:
\begin{equation}
\label{eq:rasgele2}
f(H,S_{sec})=\sup_{W \in S_{sec}} \left< H,W\right>\leq \sup_{W\in S_{sec}}\sum_{i=1}^3 \h_i^T \w_i
\end{equation}
Similarly one can achieve equality in inequality (\ref{eq:rasgele2}), although we'll not discuss here. On the other hand $W\in S_{sec}$ if and only if:
\begin{equation}
\|\w_1\|_{\ell_1}\geq\|\w_2\|_{\ell_1}
\end{equation}
We also have $\| \w_1\|_{\ell_2}^2+\| \w_2\|_{\ell_2}^2+\| \w_3\|_{\ell_2}^2=\|W\|_F^2=1$. Then we need to solve the following optimization problem (remember that $\w_i,\h_i\succeq 0~\forall~i$):
\begin{eqnarray}
\label{eq:optim2}
&&\max_{\y_1,\y_2,\y_3} ~\sum_{i=1}^3 \h_i^T \y_i\\
&&\nonumber\text{subject to}\\
&&\nonumber~~~~\y_i\succeq 0~~\forall~i\\
&&\nonumber~~~~\| \y_1\|_{\el}\geq\|\y_2\|_{\el}\\
&&\nonumber~~~~\| \y_1\|_{\ell_2}^2+\| \y_2\|_{\ell_2}^2+\| \y_3\|_{\ell_2}^2\leq1
\end{eqnarray}
Clearly, the right hand side of (\ref{eq:rasgele2}) and result of (\ref{eq:optim2}) is same again due to Lemma \ref{lem:neumann} because increasingly sorting $\y_i$'s will maximize the result. Now we'll rewrite (\ref{eq:optim2}) as follows:
\begin{eqnarray}
\label{eq:optim22}
&&\max_{\y_1,\y_2,\y_3} a_1+a_2\\
&&\nonumber\text{subject to}\\
&&\nonumber~~~~a_1=\h_1^T\y_1+\h_2^T\y_2\\
&&\nonumber~~~~a_2=\h_3^T\y_3\\
&&\nonumber~~~~\y_i\succeq 0~~\forall~i\\
&&\nonumber~~~~\| \y_1\|_{\ell_1}\geq\|\y_2\|_{\ell_1}\\
&&\nonumber~~~~\| \y_1\|_{\ell_2}^2+\|\y_2\|_{\ell_2}^2\leq E_1\\
&&\nonumber~~~~\| \y_3\|_{\ell_2}^2\leq E_2\\
&&\nonumber~~~~E_1+ E_2\leq 1
\end{eqnarray}
Now, the question is reduced to solving the following two optimization problems and maximizing over them by appropriately distributing $E_1,E_2$:
\begin{eqnarray}
\label{eq:optim23}
&&\max_{\y_1,\y_2} \h_1^T\y_1+\h_2^T\y_2\\
&&\nonumber\text{subject to}\\
&&\nonumber~~~~\| \y_1\|_{\ell_1}\geq\|\y_2\|_{\ell_1}\\
&&\nonumber~~~~\| \y_1\|_{\ell_2}^2+\|\y_2\|_{\ell_2}^2\leq 1
\end{eqnarray}
\begin{eqnarray}
\label{eq:optim24}
&&\max_{\y_3}~ \h_3^T\y_3\\
&&\nonumber\text{subject to}\\
&&\nonumber~~~~\| \y_3\|_{\ell_2}^2\leq 1
\end{eqnarray}
Let result of program (\ref{eq:optim23}) be $f_1(H,S_{sec})$ and result of program (\ref{eq:optim24}) be $f_2(H,S_{sec})$. Then clearly result of program (\ref{eq:optim22}) is
\begin{eqnarray}
\label{eq:optim25}
&&\max~ a_1+a_2\\
&&\nonumber\text{subject to}\\
&&\nonumber~~~~a_1=\sqrt{E_1}f_1(H,S_{sec})\\
&&\nonumber~~~~a_2=\sqrt{E_2}f_2(H,S_{sec})\\
&&\nonumber~~~~E_1+ E_2\leq 1
\end{eqnarray}
It is clear that $f_1(H,S_{sec}),f_2(H,S_{sec})\geq 0$. Then analyzing (\ref{eq:optim25}) we get: 
\begin{eqnarray}
\label{eq:sec1}
f(H,S_{sec})&\leq& a_1+a_2=\sqrt{E_1}f_1(H,S_{sec})+\sqrt{E_2}f_2(H,S_{sec})\\
\nonumber &\leq&\sqrt{E_1+E_2}\sqrt{f_1(H,S_{sec})^2+f_2(H,S_{sec})^2}\leq\sqrt{f_1(H,S_{sec})^2+f_2(H,S_{sec})^2}
\end{eqnarray}
Similarly one can also achieve equality in (\ref{eq:sec1}) by letting $\frac{E_1}{f_1^2}=\frac{E_2}{f_2^2}$.

Now let us turn to analyzing program (\ref{eq:optim23}). Luckily \cite{mihailo} already gives the following upper bound for this in equation (94). For any $\|\h_1\|_{\ell_1}\leq \|\h_2\|_{\ell_1}$:
\begin{equation}
\label{eq:sec2}
f_1(H,S_{sec})\leq \sqrt{\| \h_1\|_{\ell_2}^2+\| \h_2\|_{\ell_2}^2-\sum_{i=1}^ch_{2,i}^2-\frac{(\| \h_2\|_{\ell_1}-\|\h_1\|_{\ell_1}-\sum_{i=1}^c h_{2i})^2}{n-c}}
\end{equation}
for any $c\leq n-r$ such that $\| \h_2\|_{\ell_1}-\|\h_1\|_{\ell_1}-\sum_{i=1}^c h_{2,i}\geq (n-c)h_{2,c}$.

For program  (\ref{eq:optim24}) we have:
\begin{equation}
\label{eq:sec3}
f_2(H,S_{sec})\leq \h_3^T\y_3 = \left<\h_3,\y_3\right>\leq \|\h_3\|_{\ell_2}\| \y_3\|_{\ell_2}\leq \|\h_3\|_{\ell_2}
\end{equation}

Combining (\ref{eq:sec1}),(\ref{eq:sec2}),(\ref{eq:sec3}) we find:
\begin{eqnarray}
\label{eq:sec4}
f(H,S_{sec})&\leq& \sqrt{ \| H\|_F^2-\sum_{i=1}^ch_{2,i}^2-\frac{(\| \h_2\|_{\ell_1}-\|\h_1\|_{\ell_1}-\sum_{i=1}^c h_{2,i})^2}{n-c}}~~~\text{if}~\|\h_1\|_{\ell_1}\leq \|\h_2\|_{\ell_1}\\
&\leq & \|H\|_F~~~\text{else}
\end{eqnarray}

Using (\ref{eq:sec4}), for escape through a mesh (ETM) analysis we'll use the following upper bounding technique:
\begin{lem}
\label{lem:secpes1}
Let $H$ be chosen from $\G(n,n)$ and let $\h_1=\Sigma(H_{11})$ and $\h_2=\Sigma(H_{22})$. $f(H,S_{sec})=\sup_{W\in S_{sec}} \left<H,W\right>$. Then we have: $f(H,S_{sec})\leq B_{sec}$ where
\begin{align}
&\nonumber B_{sec}= \| H\|_{F}~~~\text{if }g(H,c_{sec})\leq 0\\
&\nonumber B_{sec}=\sqrt{ \| H\|_F^2-\sum_{i=1}^{c_{sec}}h_{2,i}^2-\frac{(\| \h_2\|_{\ell_1}-\|\h_1\|_{\ell_1}-\sum_{i=1}^{c_{sec}} h_{2,i})^2}{n-{c_{sec}}}}~~~\text{else}
\end{align}
where $g(H,c)=\frac{\|\h_2\|_{\ell_1}-\|\h_1\|_{\ell_1}-\sum_{i=1}^c h_{2,i}}{n-c}-h_{2,c}$ and $c_{sec}=\delta_{sec} n(1-\beta)$ is a $c\leq n(1-\beta)$ such that
\begin{align}
\label{eq:cssec}
&\nonumber c_{sec}=0~~\text{ if }\E[\|\h_2\|_{\ell_1}]\leq\E[\|\h_1\|_{\ell_1}]\\
&c_{sec}\text{ is solution of }~~(1-\eps)\frac{\E [(\|\h_2\|_{\ell_1}-\|\h_1\|_{\ell_1})-\sum_{i=1}^c h_{2,i}]}{\sqrt{n(1-\beta)}(n-c)}=F^{-1}\left(\frac{(1+\eps)c}{n(1-\beta)}\right)\text{ else }\E[\|\h_2\|_\el]>\E[\|\h_1\|_\el]
\end{align}
where $\eps>0$ can be arbitrarily small.
\end{lem}

\subsubsection{Probabilistic Analysis of $\E[B_{sec}]$}
In order to do the ETM analysis, we choose $H$ from $\G(n,n)$. Clearly $\E[B_{sec}]\geq \E[f(H,S_{sec})]$ hence we need to find an upper bound on $\E[B_{sec}]$. Similar to probabilistic analysis of strong threshold, we'll show that with high probability $g(H,c_{sec})>0$ whenever $\E[\|\h_2\|_\el-\|\h_1\|_\el]>0$. We'll declare failure else. (Failure implies $\mu=1$). Note that when $\E[\|\h_2\|_\el-\|\h_1\|_\el]>0$, we have $c_{sec}>0$.
\begin{align}
&\Prob(g(H,c_{sec})\leq0)\leq P_1+P_2~~~\text{where}\\
&P_1=\Prob(h_{2,c_{sec}}\geq \sqrt{n(1-\beta)}F^{-1}\left(\frac{(1+\eps)c_{sec}}{n(1-\beta)}\right))\\
&P_2=\Prob(\|\h_2\|_{\ell_1}-\|\h_1\|_{\ell_1}-\sum_{i=1}^{c_{sec}} h_{2,i}\leq (1-\eps)\E[\|\h_2\|_{\ell_1}-\|\h_1\|_{\ell_1}-\sum_{i=1}^{c_{sec}} h_{2,i}])
\end{align}
Remember that $h_{2,i}$ is $i$'th smallest singular value of the submatrix $H_{22}$ which is drawn from $\G(n(1-\beta),n(1-\beta))$. From quarter circle distribution it follows:
\beq
\E[h_{2,c_{sec}}]=\sqrt{n(1-\beta)}(F^{-1}(\delta_{sec})+o(1))
\eeq
Then similar to the analysis of the strong recovery using Lemmas \ref{lem:lipsc1}, \ref{lem:gauslip} and the fact that $H$ is iid Gaussian, we find
\beq
\label{eq:sec7}
P_1\leq \exp\left(-\frac{n(1-\beta)}{8}(\pi\eps\delta_{sec}+o(1))^2\right)
\eeq
Now we'll analyze $P_2$ using Lipschitzness.
\begin{lem}
\label{lem:seclip}
Let $f(H)=\|\h_2\|_{\ell_1}-\|\h_1\|_{\ell_1}-\sum_{i=1}^{c_{sec}} h_{2,i}$. Then $f$ is $\sqrt{n-c_{sec}}$ Lipschitz function of $H$.
\begin{proof}
Assume we have $2\times 2$ block matrices $\tilde{H}=\hat{H}-H\in\R^{n\times n}$ with upper left block having size $r\times r$. Also let $\|\tilde{H}\|_F=1$. Then we have 
\begin{align}
&1\geq \|\tilde{H}_{11}\|_F^2+\|\tilde{H}_{22}\|_F^2\geq \sum_{i=1}^r \sig_i(\tilde{H}_{11})^2+\sum_{i=1}^{n-r-c_{sec}}\sig_i(\tilde{H}_{22})^2\\
&\implies \sqrt{n-c_{sec}}\geq \sum_{i=1}^r \sig_i(\tilde{H}_{11})+\sum_{i=1}^{n-r-c_{sec}}\sig_i(\tilde{H}_{22})
\end{align}
Now using Lemma \ref{lem:singval} we get:
\bea
\|\tilde{H}_{11}\|_r+\|\tilde{H}_{22}\|_{n-r-c_{sec}}&\geq& \sum_{i=1}^r |\sig_i(\hat{H}_{11})-\sig_i(H_{11})|+\sum_{i=1}^{n-r-c_{sec}}|\sig_i(\hat{H}_{22})-\sig_i(H_{22})|\\
&\geq&|f(H)-f(\hat{H})|
\eea
Combining all we find: $\sqrt{n-c_{sec}}\geq |f(H)-f(\hat{H})|$ as desired.
\end{proof}
\end{lem}
Note that asymptotically $\E[f(H)]=((1-\beta)^{3/2}\gamma(1-\delta_s)-\beta^{3/2}\gamma(1)+o(1))n^{3/2}$ because $\h_1,~\h_2$ are vectors of singular values of $H_{11}$ and $H_{22}$ respectively. From (\ref{lem:seclip}) and (\ref{lem:gauslip})  we find
\beq
P_2\leq \exp\left(-\frac{n^2\eps^2}{2(1-\delta_s(1-\beta))} ((1-\beta)^{3/2}\gamma(1-\delta_s)-\beta^{3/2}\gamma(1)+o(1))^2\right)
\eeq
Finally we showed $\Prob(g(H,c_{sec})\leq0)\leq P_1+P_2$ decays to $0$ exponentially fast as $n\rightarrow\infty$. Then we use the following upper bound for $\E[B_{sec}]$:
\bea
\E[B_{sec}]\leq \int_{g(H,c_s)\leq 0} \| H\|_Fp(H)dH +\int_H \sqrt{\|H\|_F-\sum_{i=1}^{c_{sec}} h_{2,i}^2-\frac{(\|\h_2\|_{\ell_1}-\|\h_1\|_{\ell_1}-\sum_{i=1}^{c_{sec}} h_{2,i})^2}{n-c_{sec}}}p(H)dH
\eea
Using exactly same arguments in \cite{mihailo} and (\ref{eq:goeszero}) we have:
\beq
\label{eq:sec6}
\int_{g(H,c_s)\leq 0} \| H\|_Fp(H)dH \leq \exp(-\frac{n(1-\beta)}{8}(\pi\eps\delta_s+o(1))^2)n\rightarrow 0
\eeq
as $n\rightarrow \infty$. Secondly using $\sqrt{\E[\X]}\geq\E[\sqrt{\X}]$ for any R.V. $\X\geq 0$ we have:
\begin{align}
\label{eq:sec5}
&\E\left[\sqrt{\|H\|_F-\sum_{i=1}^{c_{sec}} h_{2,i}^2-\frac{(\|\h_2\|_{\ell_1}-\|\h_1\|_{\ell_1}-\sum_{i=1}^{c_{sec}} h_{2,i})^2}{n-c_{sec}}}\right]\\
&\hspace{40pt}\leq n \sqrt{1-(1-\beta)^2(1-\gamma_2(1-\delta_{sec}))-\frac{((1-\beta)^{3/2}\gamma(1-\delta_{sec})-\beta^{3/2}\gamma(1))^2}{1-\delta_{sec}(1-\beta)}+o(1)}
\end{align}
Combining (\ref{eq:sec6}), (\ref{eq:sec5}) we find that asymptotically the right hand side of (\ref{eq:sec5}) is an upper bound for $\E[B_{sec}]$. Using this we can conclude:
\begin{thm}
\label{thm:secfinal}
If $\beta\geq \frac{1}{2}$ then $\mu=1$. Otherwise:
\beq
\mu>1-(1-\beta)^2(1-\gamma_2(1-\delta_{sec}))-\frac{((1-\beta)^{3/2}\gamma(1-\delta_{sec})-\beta^{3/2}\gamma(1))^2}{1-\delta_{sec}(1-\beta)}
\eeq
is a sufficient sampling rate for $\beta$ to be sectional threshold of Gaussian operator $\A:\R^{n\times n}\rightarrow \R^{\mu n^2}$, where $0<\delta_{sec}<1$ is solution of (from (\ref{eq:cssec}))
\beq
(1-\eps)\frac{(1-\beta)^{3/2}\gamma(1-\delta)-\beta^{3/2}\gamma(1)}{\sqrt{1-\beta}(1-\delta(1-\beta))}=F^{-1}((1+\eps)\delta)
\eeq
\end{thm}
In order to find the least $\mu$ we let $\eps\rightarrow 0$. Numerical calculations result in sectional threshold of Figure (\ref{fig:real_thresholds}).

\subsection{Weak Threshold}

In this section, we'll derive the relation between $\mu$ and $\beta$ for the weak threshold described below.
\begin{WeakDef}
Let $\A:\R^{n\times n}\rightarrow \R^m$ be a random Gaussian operator and let $X\in\R^{n\times n}$ be an arbitrary matrix with $\text{rank}(X)=\beta n$. We say that $\beta$ is a weak recovery threshold if with high probability this particular matrix $X$ can be recovered from measurements $\A(X)$ via program (\ref{eq: nuc norm min}).
\end{WeakDef}

We remark that the weak threshold is the one that can be observed from simulations. The strong (and sectional) thresholds cannot because there is no way to check the recovery of {\bf{all}} low rank $X$ (or all $X$ of a particular support). In this sense, the weak threshold is the most important.

Again given a fixed $\beta$, we'll aim to the least $\mu$ such that $\beta$ is weak threshold for a random Gaussian operator $\A:\R^{n\times n}\rightarrow \R^{\mu n^2}$. In order to prevent repetitions, we'll be more concise in this section because many of the derivations are repetitions of the derivations for strong and sectional thresholds.

\begin{lem}
\label{lem:weak_null}
Let $X\in\R^{n\times n}$ with $\text{rank}(X)=r$, SVD $X=U\Sigma V^T$ with $\Sigma\in\R^{r\times r}$. Then it can be recovered using (\ref{eq: nuc norm min}) iff for all $W\in \N(\A)$ we have
\begin{equation}
\label{eq: weak null}
\text{trace}(U^TWV)+\|\bar{U}^TW\bar{V}\|_\s>0
\end{equation}
where $\bar{U},\bar{V}$ such that $[U~\bar{U}]$ and $[V~\bar{V}]$ are unitary.
\begin{proof}
Since singular values are unitarily invariant , if (\ref{eq: weak null}) holds using Lemma (\ref{lem:simple}):
\bea
\| X+W\|_\s &=&\| [U~\bar{U}]^T (X+W) [V~\bar{V}]\|_\s=\left\|\begin{bmatrix} \Sigma & 0\\0 & 0\end{bmatrix}+\begin{bmatrix} U^TWV & \dots\\\dots & \bar{U}^TW\bar{V}\end{bmatrix}\right\|_\s\\
&\geq&\| \Sigma+U^TWV\|_\s+\|\bar{U}^TW\bar{V}\|_\s\geq \text{trace}(\Sigma+U^TWV)+\|\bar{U}^TW\bar{V}\|_\s\\
&\geq& \| X\|_\s+\text{trace}(U^TWV)+\|\bar{U}^TW\bar{V}\|_\s>\| X\|_\s
\eea
Hence $X$ is unique minimizer of program (\ref{eq: nuc norm min}). In \cite{arxiv}, it was shown that condition (\ref{eq: weak null}) is tight in the sense that if there is a $W\in\N(\A)$ such that $\text{trace}(U^TWV)+\|\bar{U}^TW\bar{V}\|_\s<0$ then $X$ is not minimizer of (\ref{eq: nuc norm min}).
\end{proof}
Note that conditions (\ref{eq: weak null}) is independent of singular values of $X$. This suggests that not only $X$ but also all matrices with same left and right singular vectors $U,V$ are recoverable via (\ref{eq: nuc norm min}).
\end{lem}

{\bf{Analyzing the condition:}} $\A$ is a random Gaussian operator and we'll analyze the linear regime where $m=\mu n^2$ and $r=\beta n$ with $\beta>0$.

Let $S_w$ be the set of all matrices such that $\text{trace}(U^TWV)+\|\bar{U}^TW\bar{V}\|_\s\leq0$ and $\| W\|_F=1$. We need to make sure null space $\N(\A)$ has no intersection with $S_w$. We first upper bound (\ref{eq: ETMeq}) in Theorem \ref{thm: ETM}, then choose $m$ (and $\mu$) respectively. As discussed in \cite{arxiv}, without loss of generality we can assume $X=\begin{bmatrix} \Sigma & 0\\0&0\end{bmatrix}$ where $\Sigma\in\R^{r\times r}$ is diagonal matrix with positive diagonal. This is because any $X=U\Sigma V^T$ can be transformed into this form by unitary transformation $[U~\bar{U}]^T X[V~\bar{V}]$ and since null space is uniformly chosen (i.e. its basis is $n^2-m$ random matrices iid chosen from $\G(n,n)$) after unitary transformation its distribution will still be uniform. Then $S_w$ can be assumed to be set of matrices with $\text{trace}(W_{11})+\|W_{22}\|_\s\leq 0$ and $\|W\|_F=1$.

Again we assume $2\times 2$ block matrices. Firstly, given a fixed $H\in \R^{n\times n}$ we'll calculate an upper bound on $f(H,S_w)=\sup_{W_\in S_w} \left< H,W\right> $.

Let $\h_1=diag(H_{11})$ i.e. diagonal entries of $H_{11}$, $\h_2=\Sigma(H_{22})$ and let $\h_3$ be increasingly sorted absolute values of remaining entries, which are entries of $H_{12},H_{21}$ and off-diagonal entries of $H_{11}$. $\w_1,\w_2,\w_3$ is defined similarly. Also let $x_{i,j}$ denote $j$'th entry of vector $\x_i$.

From Lemma (\ref{lem:neumann}) we have:
\begin{equation}
\label{eq:rasgele3}
f(H,S_w)=\sup_{W \in S_w} \left< H,W\right>\leq \sup_{W\in S_w}\sum_{i=1}^3 \h_i^T \w_i
\end{equation}

We introduce the following notation. Let $s(\x)$ denote summation of entries of $\x$ i.e. $s(\x)=\sum_i x_i$. Then $W\in S_w$ if and only if:
\begin{align}
&s(\w_1)+\|\w_2\|_{\ell_1}\leq 0~~\text{and}\\
&\| \w_1\|_{\ell_2}^2+\|\w_2\|_{\ell_2}^2+\|\w_3\|^2_{\ell_2}=\|W\|_F^2= 1
\end{align}
Then we need to solve the following equivalent optimization problem given $H$  (Note that $\w_i,\h_i\succeq 0~\forall~2\leq i\leq 3$):
\begin{eqnarray}
\label{eq:optim3}
&&\max_{\y_1,\y_2,\y_3} \sum_{i=1}^3 \h_i^T \y_i\\
&&\nonumber\text{subject to}\\
&&\nonumber~~~~\y_2,\y_3\succeq0\\
&&\nonumber~~~~s(\y_1)+\|\y_2\|_{\ell_1}\leq 0\\
&&\nonumber~~~~\|\y_1\|_{\ell_2}^2+\|\y_2\|_{\ell_2}^2+\|\y_3\|_{\ell_2}^2\leq 1
\end{eqnarray}
Right hand side of (\ref{eq:rasgele3}) and output of (\ref{eq:optim3}) is same. We'll rewrite (\ref{eq:optim3}) as follows:
\begin{eqnarray}
\label{eq:optim32}
&&\max_{\y_1,\y_2,\y_3} a_1+a_2\\
&&\nonumber\text{subject to}\\
&&\nonumber~~~~a_1=\h_1^T\y_1+\h_2^T\y_2\\
&&\nonumber~~~~a_2=\h_3^T\y_3\\
&&\nonumber~~~~\y_2,\y_3\succeq0\\
&&\nonumber~~~~s(\y_1)+\|\y_2\|_{\ell_1}\leq 0\\
&&\nonumber~~~~\| \y_1\|_{\ell_2}^2+\|\y_2\|_{\ell_2}^2\leq E_1\\
&&\nonumber~~~~\| \y_3\|_{\ell_2}^2\leq E_2\\
&&\nonumber~~~~E_1+ E_2\leq 1
\end{eqnarray}
Note that (\ref{eq:optim32}) is essentially same as (67) of \cite{mihailo}. Basically (\ref{eq:optim32}) has additional terms of $\h_3,\y_3$ and $y_{1,i}$ corresponds $-y_{n-r+i}$ of \cite{mihailo} for $1\leq i\leq r$ and $y_{2,j}$ corresponds $y_j$ of \cite{mihailo} for $1\leq j\leq n-r$.  Then repeating exactly same steps that come before Lemma (\ref{lem:secpes1}) and equation (\ref{eq:sec4}) and using (67), (68) of \cite{mihailo} we find:
\begin{lem}
\label{lem:weak1}
\begin{eqnarray}
f(H,S_w)&\leq& \sqrt{ \| H\|_F^2-\sum_{i=1}^ch_{2,i}^2-\frac{(\| \h_2\|_{\ell_1}+s(\h_1)-\sum_{i=1}^c h_{2,i})^2}{n-c}}~~~\text{if}~s(\h_1)+\|\h_2\|_\el> 0\\
&\leq & \|H\|_F~~~\text{else}
\end{eqnarray}
for any $0\leq c\leq n-r$ such that $s(\h_1)+\|\h_2\|_\el-\sum_{i=1}^c h_{2,i}\geq (n-c)h_{2,c}$.
\end{lem}

Based on Lemma \ref{lem:weak1}, for ETM analysis we'll use the following lemma:
\begin{lem}
\label{lem:weak2}
Let $H$ be chosen from $\G(n,n)$ and let $\h_1=diag(H_{11})$, $\h_2=\Sigma(H_{22})$ and $f(H,S_w)=\sup_{W\in S_w} \left<H,W\right>$. Then we have: $f(H,S_w)\leq B_w$ where
\begin{align}
&\nonumber B_w= \| H\|_{F}~~~\text{if }g(H,c_w)\leq 0\\
&\nonumber B_w=\sqrt{ \| H\|_F^2-\sum_{i=1}^{c_w}h_{2,i}^2-\frac{(\| \h_2\|_\el+s(\h_1)-\sum_{i=1}^{c_w} h_{2,i})^2}{n-{c_w}}}~~~\text{else}
\end{align}
where $g(H,c)=\frac{\|\h_2\|_{\ell_1}+s(\h_1)-\sum_{i=1}^c h_{2,i}}{n-c}-h_{2,c}$ and $c_w=\delta_w n(1-\beta)$ is a $c\leq n(1-\beta)$ such that
\begin{align}
\label{eq:csweak}
\nonumber c_w&=0~~~\text{ if }\E[\|\h_2\|_\el+s(\h_1)]\leq0\\
c_w&\text{ is solution of }~~~(1-\eps)\frac{\E [\|\h_2\|_\el+s(\h_1)-\sum_{i=1}^c h_{2,i}]}{\sqrt{n(1-\beta)}(n-c)}=F^{-1}\left(\frac{(1+\eps)c}{n(1-\beta)}\right)\text{ else }\E[\|\h_2\|_{\ell_1}+s(\h_1)]>0
\end{align}
where $\eps>0$ can be arbitrarily small.
\end{lem}
Note that, when $H$ is iid Gaussian, for any $\beta<1$ we have $\E[\|\h_2\|_{\ell_1}+s(\h_1)]>0$ since entries of $\h_1$ is iid Gaussian hence $\E[s(\h_1)]=0$ and clearly $\E[\|\h_2\|_\el]>0$. As a result $c_w>0$ too.

\subsubsection{Probabilistic Analysis of $\E[B_w]$}
\begin{figure}[t]
\centering
  \includegraphics[width= 1\textwidth]{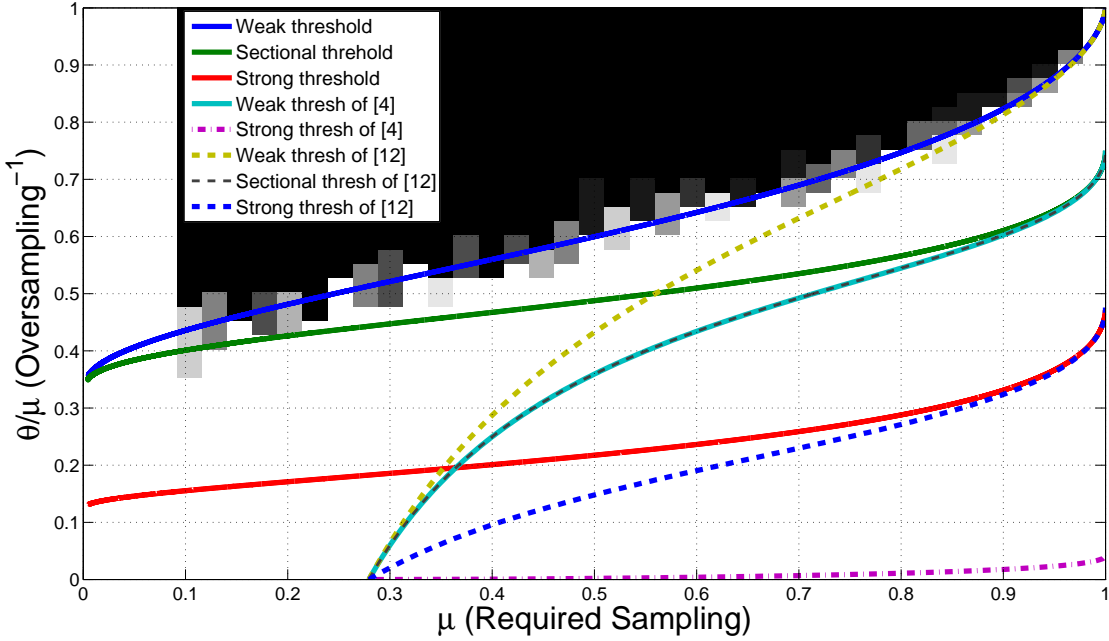}
 \caption{\scriptsize{Results of \cite{Recht_Xu_Hassibi,arxiv} vs results we get by using escape through a mesh (ETM) analysis similar to \cite{mihailo}} (For square matrices). Here $\theta$ is model complexity i.e. degrees of freedom of the matrix, ($\theta=\beta(2-\beta)$). This plot gives the efficiency of nuclear norm minimization as a function of number of samples $\mu$. It gives at a certain $\mu$, how much more one should oversample the content of the matrix to perform NNM successfully. Simulations are done for $40\times 40$ matrices and program (\ref{eq: nuc norm min}) is solved with Gaussian measurements. Our weak threshold and simulations match almost exactly. Black regions indicate failure and white regions mean success. Due to low precision ($40$ is small), we did not include simulation results for $\mu\leq 0.1$.}
  \label{fig:real_thresholds}
\end{figure}

Similar to previous analysis $H$ is drawn from $\G(n,n)$ and in order to use Theorem \ref{thm: ETM} we need to upper bound $\E[B_w]$. Using the same steps and letting $f(H)=\|\h_2\|_\el+s(\h_1)-\sum_{i=1}^{c_w} h_{2,i}$:
\begin{align}
&\Prob(g(H,c_w)\leq0)\leq P_1+P_2~~~\text{where}\\
&P_1=\Prob(h_{2,c_w}\geq \sqrt{n(1-\beta)}F^{-1}\left(\frac{(1+\eps)c_w}{n(1-\beta)}\right))\\
&P_2=\Prob(f(H)\leq (1-\eps)\E[f(H)])
\end{align}
An upper bound for $P_1$ was already given in (\ref{eq:sec7}). Also similar to Lemma (\ref{lem:seclip}) one can show $f(H)$ is $\sqrt{n-c_w}$ Lipschitz function of $H$. Therefore $g(H,c_w)$ will approach to $0$ exponentially fast ($e^{-O(n)}$). Combining this and the same arguments prior to (\ref{eq:sec5}) yields:
\begin{eqnarray}
\E[B_w]&\leq& \sqrt{\E\left[\| H\|_F^2-\sum_{i=1}^{c_w}h_{2,i}^2-\frac{(\| \h_2\|_\el+s(\h_1)-\sum_{i=1}^{c_w} h_{2,i})^2}{n-{c_w}}\right]}+o(1)\\
&\leq&n\sqrt{1-(1-\beta)^2(1-\gamma_2(1-\delta_w))-\frac{(1-\beta)^3\gamma(1-\delta_w)^2}{1-\delta_w(1-\beta)}}+o(1)
\end{eqnarray}
Then using Theorem \ref{thm: ETM} and (\ref{eq:csweak}) we can write:
\begin{thm}
\label{thm:weakfinal}
\beq
\mu>1-(1-\beta)^2(1-\gamma_2(1-\delta_w))-\frac{(1-\beta)^3\gamma(1-\delta_w)^2}{1-\delta_w(1-\beta)}
\eeq
is a sufficient sampling rate for $\beta$ to be weak threshold of Gaussian operator $\A:\R^{n\times n}\rightarrow\R^{\mu n^2}$ where $0<\delta_w<1$ is solution of
\beq
(1-\eps)\frac{(1-\beta)^{3/2}\gamma(1-\delta)}{\sqrt{1-\beta}(1-\delta(1-\beta))}=F^{-1}((1+\eps)\delta)
\eeq  
\end{thm}
To find the least $\mu$ we let $\eps\rightarrow 0$. Numerical calculations result in weak threshold of Figure (\ref{fig:real_thresholds}).

\section{Thresholds for Positive Semidefinite Matrices}

\subsection{Additional Notations and Lemmas}
Before starting our analysis, we'll briefly introduce some more notations and lemmas.

$\PS_n$ denotes the set of Hermitian (real and symmetric) matrices of size $n\times n$. Similarly $\PS_+^n$ denotes the set of positive semidefinite matrices. PSD stands for positive semidefinite.

Let $\N_s(\A)\subset\N(\A)$ denote subspace of null space of $\A$ which consists of Hermitian matrices. 

Denote Gaussian unitary ensemble by $\HH(n)$ which is ensemble of Hermitian matrices of size $n \times n$ with independent Gaussian entries in the lower triangular part, where off-diagonal entries have variance $1$ and diagonal entries have variance $2$. In order to create such a matrix $B$ one can choose a matrix $A$ from $\G(n,n)$ and then let $B=\frac{A+A^T}{\sqrt{2}}$.

Let $X\in\PS^n$ with $\text{rank}(X)=r$. Then (skinny) eigenvalue decomposition (EVD) of $X$ is $X=U\Lambda U^T$ for some partial unitary $U\in\R^{n\times r}$ and diagonal matrix $\Lambda\in\R^{r\times r}$. Denote $i$'th largest eigenvalue of $X$ by $\lambda_i(X)$ for $1\leq i\leq n$. Let $\Lambda(X)$ denote increasingly ordered eigenvalues of $X\in\PS^n$ i.e. $\Lambda(X)=[\lambda_n(X)~\dots~\lambda_1(X)]$. Also note that singular values of $X$ corresponds to absolute values of eigenvalues of $X$.

Let $c=\sqrt{1/2}$. If $A\in \PS^n$ is a symmetric real matrix, define $vec(A)\in \R^{n(n+1)/2}$ to be following vector:
\beq
vec(A)=\frac{1}{c}[cA_{1,1}~A_{2,1}~\dots~A_{n,1}~cA_{2,2}~A_{3,2}~\dots~A_{n,2}~cA_{3,3}~A_{3,4}~\dots~cA_{n-1,n-1}~A_{n,n-1}~cA_{nn}]^T
\eeq
In other words for each $i$ we let $b_i=[cA_{i,i}~A_{i+1,i}\dots~A_{n,i}]/c$ and we let $vec(A)=[b_1~b_2\dots~b_n]^T$.

Now note that for any $A,B\in\PS^n$ we have
\beq
\left<A,B\right>=\sum_{i,j} A_{i,j}B_{i,j}=\sum_i A_{i,i}B_{i,i}+2\sum_{i< j} A_{i,j}B_{i,j}=vec(A)^Tvec(B)
\eeq
Clearly $vec(.):\PS^n\rightarrow \R^{n(n+1)/2}$ is bijective. Then let $ivec(.)$ denote inverse of the function $vec(.)$. Also it is clear that $vec(.)$ is linear.

Let $\eta_-(X),\eta_+(X),\eta_0(X)$ denote number of positive, negative and zero eigenvalues of $X$. The triple $(\eta_-(X),\eta_0(X),\eta_+(X))$ is called the inertia of $X$.

Similar to $\G(n,n)$ following limits exist for a random matrix $X$ drawn from $\HH(n)$:
The histogram of the eigenvalues of $X$ normalized by $1/\sqrt{n}$ converges to semicircle distribution given by:
\bea
\phi_s(x)&=&\frac{\sqrt{4-x^2}}{2\pi}~~~\text{for}~-2\leq x\leq 2\\
\phi_s(x)&=&0~~~\text{else}.
\eea
For $-2\leq x\leq 2$, define the cumulative distribution function corresponding to the semicircle law as:
\beq
F_s(x)=\int_{-2}^x\phi_s(t)dt
\eeq
Similar to the definitions of $\gamma(.)$ and $\gamma_2(.)$ we can define $\gamma_s(\beta)$ and $\gamma_{2,s}(\beta)$ for $0\leq \beta\leq 1$: 
\begin{align}
&\gamma_s(\beta)=\lim_{n\rightarrow \infty}\frac{\sum_{i=1}^{\beta n} \la_i(X)}{n^{3/2}}\\
&\gamma_{2,s}(\beta)=lim_{n\rightarrow \infty}\frac{\sum_{i=1}^{\beta n} \la_i(X)^2}{n^2}
\end{align}
Using definitions of $F(.)$ and $F_s$ it is easy to see that
\begin{align}
&\gamma_s(\beta)=\frac{\gamma(2\beta)}{2}~~~\text{for}~\beta\leq 0.5\\
&\gamma_s(\beta)=\frac{\gamma(2-2\beta)}{2}~~~\text{for}~0.5\leq \beta\leq 1
\end{align}
Similarly we have
\begin{align}
&\gamma_{2,s}(\beta)=\frac{\gamma_2(2\beta)}{2}~~~\text{for}~\beta\leq 0.5\\
&\gamma_{2,s}(\beta)=1-\frac{\gamma_2(2-2\beta)}{2}~~~\text{for}~0.5\leq \beta\leq 1
\end{align}

We'll need the following lemmas about eigenvalues of matrices:
\begin{lem}
\label{lem:eiglem}
Let $X,Y\in \PS^n$, $s_i(X,Y)=|\la_i(X)-\la_i(Y)|$ and let $s_{[1]}(X,Y)\geq s_{[2]}(X,Y)\geq \dots \geq s_{[n]}(X,Y)$ be a decreasingly ordered arrangement of $\{s_i(X,Y)\}_{i=1}^n$. Then from \cite{lidskii} we have the following inequality:
\beq
\sum_{i=1}^ks_{[i]}(X,Y)\leq \sum_{i=1}^k \sig_i(X-Y)
\eeq
for any $1\leq k\leq n$. This is the eigenvalue counterpart of (\ref{lem:singval}).
\end{lem}
\begin{lem}
\label{lem:eiglip}
$\la_i(X)$ is $1$-Lipschitz function of $X$ (and $vec(X)$).
\begin{proof}
Given $X,Y$ we have:
\beq
\|X-Y\|_F\geq \sig_1(X-Y)\geq s_{[1]}(X,Y)\geq s_i(X,Y)=|\la_i(X)-\la_i(Y)|
\eeq
\end{proof}
\end{lem}

\begin{lem}
\label{lem:psdsimple}
If $X,Y\in\R^{n\times n}$ are positive semidefinite matrices then $\text{trace}(XY)\geq 0$.
\begin{proof}
Let $A,B\in\R^{n\times n}$ be arbitrary square-roots of $X,Y$ respectively. In other words $A^TA=X$ and $B^TB=Y$. Since $X,Y$ are PSD $A,B$ exists. Then we can write:
\beq
\text{trace}(XY)=\text{trace}(A^TAB^TB)=\text{trace}(AB^TBA^T)=\text{trace}(AB^T(AB^T)^T)=\|AB^T\|_F^2\geq 0
\eeq
\end{proof}
\end{lem}

\begin{lem}
\label{lem:inertia}
Given $A,B\in\PS^n$ we have following inequalities due to \cite{greg}
\bea
&\eta_+(A)-\eta_-(B)\leq \eta_+(A+B)\leq \eta_+(A)+\eta_+(B)\\
&\eta_-(A)-\eta_+(B)\leq \eta_-(A+B)\leq \eta_-(A)+\eta_-(B)
\eea
\end{lem}

\subsection{PSD Recovery Methods}
Now we'll state and analyze null space conditions for success of the following program which is equivalent to nuclear norm minimization for PSD matrices. $X_0\in\PS_+^n$ be a (low rank) matrix. Then we want $X_0$ to be unique solution of following program:
\begin{eqnarray}
\label{eq: trace min}
&&\min~\text{trace}(X)\\
&&\mbox{subject to}\\
 &&\hspace{20 pt}\A(X)=\A(X_0), \nonumber\\
&&\hspace{20 pt}\nonumber X\succeq 0
\end{eqnarray}
This is equivalent to (\ref{eq: nuc norm min}) because $\text{trace}(X_0)=\sum_{i=1}^n\la_i(X)=\|X_0\|_\s$ since eigenvalues and singulars are same for PSD matrices. Similar to previous discussion, measurement operator is random Gaussian.

In addition to this, we'll state the results for the following program where we want $X_0$ to be unique positive semidefinite solution satisfying measurements $A(X_0)$:
\begin{eqnarray}
\label{eq: unique PSD}
&&\text{find}~X\\
&&\mbox{subject to}\\
 &&\hspace{20 pt}\A(X)=\A(X_0), \nonumber\\
&&\hspace{20 pt}\nonumber X\succeq 0
\end{eqnarray}
However we'll omit the analysis for this, because it is very similar to program (\ref{eq: trace min}) and actually simpler to analyze.

\subsection{PSD Weak Threshold}
\begin{PSDWeak}
$\beta$ is called a PSD weak threshold for random Gaussian operator $\A:\R^{n\times n}\rightarrow \R^{\mu n(n+1)/2}$ if given a fixed $X\in\PS_+^n$ with $\text{rank}(X)=\beta n$, $X$ can be recovered from measurements $\A(X)$ via (\ref{eq: trace min}) asymptotically with probability $1$.

\end{PSDWeak}
For a given $\beta$, our aim is to find the least $\mu\leq 1$ so that $\beta$ is a weak threshold for $\A$.

\begin{lem}
\label{lem:PSDweak}
Let $X\in\PS_+^n$ be a rank $r$ matrix with EVD $X=U\Lambda U^T$ and $\Lambda\in\R^{r\times r}$. Then $X$ is unique minimizer of (\ref{eq: trace min}) if for all $W\in\N(\A)$ we have:
\begin{align}
\label{eq:propw}
&W~\text{is not hermitian}~~~\text{or}\\
\label{eq:propw2}
&\text{trace}(W)>0~~~\text{or}\\
\label{eq:propw3}
&\eta_-(\bar{U}^TW\bar{U})>0
\end{align}
\end{lem}
where $[U~\bar{U}]$ is a unitary matrix.
\begin{proof}
If $W$ is not hermitian then $X+W$ is not hermitian thus not PSD. If $\text{trace}(W)>0$ then $\text{trace}(X+W)>\text{trace}(X)$ as desired. On the other hand if $\bar{U}^TW\bar{U}$ has a negative eigenvalue, we can write
\beq
Y:=[U~\bar{U}]^T (X+W) [U~\bar{U}]=\begin{bmatrix} \Lambda+U^TWU & \dots\\\dots&\bar{U}^TW\bar{U}\end{bmatrix}
\eeq
which means lower right submatrix $\bar{U}^TW\bar{U}$ of $Y$ (which is a principal submatrix) is not PSD. Then it immediately follows that $Y$ is not PSD, because we can find a vector $\v\in\R^n$ to make $\v^TY\v<0$. Then $X+W$ is not PSD as it can be obtained by unitarily transforming $Y$ (i.e. $[U~\bar{U}]Y[U~\bar{U}]^T$) which preserves eigenvalues. Then as long as $W$ satisfies one of the (\ref{eq:propw}), $X+W$ cannot be minimizer hence $X$ is unique minimizer.
\end{proof}
One can also give the if and only if condition for PSD weak recovery. Without proof, we'll state the difference from Lemma \ref{lem:PSDweak}:

For all $W\in\N(\A)$ we should have: (\ref{eq:propw}) or (\ref{eq:propw2}) or (\ref{eq:propw3}) or {\bf{"column space of $\bar{U}^TWU$ is not a subset of column space of $\bar{U}^TW\bar{U}$"}}.

However, this last condition (in bold) would not have any affect in our ETM analysis. (Again without proof) The reason is that, with arbitrarily small perturbation we can make $\bar{U}^TW\bar{U}$ full rank, while not changing other properties of $W$ at all, hence the last condition will be obsolete.

\begin{lem}
Conditions (\ref{eq:propw}, \ref{eq:propw2}, \ref{eq:propw3}) is also sufficient to guarantee sectional recovery. In other words given an $X=U\Lambda U^T$, if (\ref{eq:propw}, \ref{eq:propw2}, \ref{eq:propw3}) holds for all $W\in\N(\A)$, then in addition to recoverability of $X$, we can recover all PSD matrices $Y$ with support $\{UU^T,UU^T\}$ from measurements $\A(Y)$ with (\ref{eq: trace min}). 
\begin{proof}
Any PSD $Y$ with support $UU^T$ can be written as $Y=U_Y\Lambda_YU_Y^T$ with $U_YU_Y^T=UU^T$. Now, assume (\ref{eq:propw}, \ref{eq:propw2}, \ref{eq:propw3}) holds. Then, if we have $\eta_-(\bar{U}_Y^TW\bar{U}_Y)>0$ whenever $\eta_-(\bar{U}^TW\bar{U})>0$, using Lemma (\ref{lem:PSDweak}) we are done, because all conditions for $Y$ become satisfied.

As a result, it remains to show: $\eta_-(\bar{U}_Y^TW\bar{U}_Y)>0\iff \eta_-(\bar{U}^TW\bar{U})>0$
\begin{proof}
Let $\v\in\R^{n-r}$ such that $\v^T\bar{U}^TW\bar{U}\v<0$. Then since column spaces of $\bar{U}$ and $\bar{U}_Y$ are same we can choose $\v_2=\bar{U}_Y^T\bar{U}\v$ so that $\bar{U}_Y\v_2=\bar{U}\v\implies\v_2^T\bar{U}_Y^TW\bar{U}_Y\v_2<0$.
\end{proof}\end{proof}
\end{lem}
This result suggests that there is no need to analyze sectional condition separately because results we get for weak will also work for sectional.\\

 Now we'll start null space analysis. Let $\A:\R^{n\times n}\rightarrow \R^{m}$ be random Gaussian operator where $m=\mu n(n+1)/2$ ($0\leq \mu\leq 1$). In \cite{weiyu} it was argued that distribution of $\N_s(\A)$ (null space restricted to hermitians) is equivalent to a subspace having matrices $\{D_i\}_{i=1}^{n(n+1)/2-m}$ as basis where $\{D_i\}_{i=1}^{n(n+1)/2-m}$ is drawn iid from $\HH(n)$. This is easy to see when we consider $\A(\cdot)$ as a mapping from lower triangular entries to $\R^m$.
 
This also implies distribution of $\N_s(\A)$ is unitarily invariant because if $D$ is chosen from $\HH(n)$ then for a fixed unitary matrix $V$, $VDV^T$ and $D$ has same distribution (identical random variables).
\begin{proof}
$D$ is equivalent to $(G+G^T)/\sqrt{2}$ where $G$ is chosen from $\G(n,n)$. Then $VDV^T$ is equivalent to $(VGV^T+(VGV^T)^T)/\sqrt{2}$. Now using distribution of $VGV^T$ is equivalent to that of $G$ we end up with the desired result.
\end{proof}

Let $X=U\Lambda U^T$ be given where $\text{rank}(X)=r=\beta n$. Similar to previous analysis let $S_{wp}$ denote the set of hermitian matrices $W$ so that $\text{trace}(W)\leq 0$, $\eta_-(\bar{U}^TW\bar{U})=0$ and $\|W\|_F=1$. Since $\N_s(\A)$ is unitarily invariant we can assume $X$ is diagonal and: $X=\begin{bmatrix} \Lambda & 0\\ 0 & 0\end{bmatrix}$.
 Now condition $\eta_-(\bar{U}^TW\bar{U})=0$ can be replaced by $\eta_-(W_{22})=0$. We want to make sure that $\N_s(\A)$ does not intersect with $S_{wp}$ so that null space condition (\ref{eq:propw}) will be satisfied.

Note that this can be rewritten in the following way which will enable us to use Theorem \ref{thm: ETM}. Let $vec(\N_s(\A))$ be the subspace formed by applying $vec(.)$ to elements of $\N_s(\A)$. Then its distribution is equivalent to a subspace in $\R^{n(n+1)/2}$ having basis $\{\dd_i\}_{i=1}^{(1-\mu)n(n+1)/2}$ where $\{\dd_i\}_i$'s are iid vectors drawn from $\G(n(n+1)/2,1)$. This is because random vectors $\{vec(D_i)\}_i$ are equivalent to $\{\dd_i\}_i$ scaled by $\sqrt{2}$. Similarly let $vec(S_{wp})$ be the set of vectors obtained by applying $vec(.)$ to elements of $S_{wp}$.

Then because $vec(.)$ is linear, $\N_s(\A)\cap S_{wp}= \emptyset \iff vec(\N_s(\A))\cap vec(S_{wp})= \emptyset$. Let $\h$ be drawn from $\G(n(n+1)/2,1)$ and let $D$ be drawn from $\HH(n)$. Since $vec(.)$ also preserves inner products, in order to use Theorem (\ref{thm: ETM}) as previously we just need to calculate:
\beq
\omega(S_{wp})=\E[\sup_{\w\in vec(S_{wp})} \h^T\w]=\frac{1}{\sqrt{2}}\E[\sup_{W\in S_{wp}} \left<D,W\right>]
\eeq

Let $H$ be Hermitian and define $f(H,S_{wp})=\sup_{W\in S_{wp}} \left<H,W\right>$. Then we'll firstly upper bound $f(H,S_{sw})$ then take expectation of upper bound as previously. Let $s(\x)$ denote summation of entries of vector $\x$.

Let $\h_1$ denote the diagonal entries of $H_{11}$, $\h_2=\Lambda(H_{2,2})$ and $\h_3$ denote increasingly ordered absoute values of entries of $H_{12},H_{21}$ and off-diagonal entries of $H_{11}$. $\w_1,\w_2,\w_3$ are defined similarly for $W$. Note that $W\in S_{sw}$ if and only if
\begin{align}
&\w_2\succeq 0~~~\text{(i.e.}~W_{2,2}~\text{PSD)}\\
&s(\w_1)+s(\w_2)\leq 0\\
&\|\w_1\|_{\ell_2}^2+\|\w_2\|_{\ell_2}^2+\|\w_3\|_{\ell_2}^2=1
\end{align}
Now using Lemma (\ref{lem:neumann}) we write:
\beq
\left<H,W\right>\leq \h_1^T\w_1+\h_3^T\w_3+\left<H_{22},W_{22}\right>
\eeq
Let $H_{2,2}=H^+_{2,2}-H^-_{2,2}$ where both of $H^+_{2,2},H^-_{2,2}\succeq 0$. Then from Lemma \ref{lem:psdsimple} we get $\text{trace}(W_{2,2}^TH^-_{2,2})\geq 0$ and from Lemma (\ref{lem:neumann}) we get: $\text{trace}(W_{2,2}^TH^+_{2,2})\leq \sum_{i=1}^{\eta_+(H_{2,2})}\la_i(W_{2,2})\la_i(H^+_{2,2})$. Combining these we can find:
\beq
\left<H_{22},W_{22}\right>\leq \sum_{i=1}^{\eta_+(H_{2,2})}\la_i(W_{2,2})\la_i(H^+_{2,2})
\eeq
Clearly $\la_i(H^+_{2,2})=\la_i(H_{2,2})$ for any $i\leq \eta_+(H_{2,2})$. Now let $\h_4,\w_4$ be increasingly ordered first $\eta_+(H_{2,2})$ eigenvalues of $H_{2,2},W_{2,2}$ respectively. Then clearly $\|\w_4\|_{\ell_2}\leq \|\w_2\|_{\ell_2}$ and $s(\w_4)\leq s(\w_2)$ since $W_{2,2}$ is PSD. Then an upper bound for $\left<H,W\right>$ is $\h_1^T\w_1+\h_3^T\w_3+\h_4^T\w_4$.

Remember that $\h_3,\h_4,\w_2,\w_3,\w_4\succeq 0$ in the previous discussions. Then solution of the following program will give the upper bound for $f(H,S_{wp})$:
\begin{eqnarray}
\label{eq:optim4}
&&\max_{\y_1,\y_3,\y_4} \h_1^T\y_1+\h_3^T\y_3+\h_4^T\y_4\\
&&\nonumber\text{subject to}\\
&&\nonumber~~~~\y_3,\y_4\succeq0\\
&&\nonumber~~~~s(\y_1)+s(\y_4)\leq 0\\
&&\nonumber~~~~\|\y_1\|_{\ell_2}^2+\|\y_3\|_{\ell_2}^2+\|\y_4\|_{\ell_2}^2\leq 1
\end{eqnarray}
Similar to this, following will also yield the same upper bound:
\begin{eqnarray}
\label{eq:optim41}
&&\max_{\y_1,\y_2,\y_3} \h_1^T\y_1+\h_2^T\y_2+\h_3^T\y_3\\
&&\nonumber\text{subject to}\\
&&\nonumber~~~~\y_2,\y_3\succeq0\\
&&\nonumber~~~~s(\y_1)+s(\y_2)\leq 0\\
&&\nonumber~~~~\|\y_1\|_{\ell_2}^2+\|\y_2\|_{\ell_2}^2+\|\y_3\|_{\ell_2}^2\leq 1
\end{eqnarray}
The difference is $\h_2\not\succeq 0$ however largest $\eta_+(H_{2,2})$ entries of $\h_2$ (i.e. positive entries) gives $\h_4$ and the remaining entries are nonpositive. Then programs (\ref{eq:optim4}) and (\ref{eq:optim41}) gives the same result since in order to maximize $\h_2^T\y_2$ one should set $y_{2,i}=0$ whenever $h_{2,i}\leq 0$ since we need to satisfy $\y_2\succeq 0$. In other words for any $\y_2\succeq 0$, the vector $\v=[y_{2,1}~\dots~y_{2,\eta_+(H_{2,2})}~0~0~\dots~0]^T$ is feasible and it will yield better or equal result because we'll have: $\h_2^T\v\geq \h_2^T\y_2$, $s(\v)\leq s(\y_2)$ and $\|\v\|_{\ell_2}\leq\|\y_2\|_{\ell_2}$. Consequently (\ref{eq:optim41}) reduces to (\ref{eq:optim4}).

Now note that (\ref{eq:optim1}) is exactly same as program (115) of \cite{mihailo} except additional terms of $\h_3,\y_3$. Then using (116) of \cite{mihailo} and repeating the steps before (\ref{eq:sec4}) and using $\|\h_1\|_{\ell_2}^2+\|\h_2\|_{\ell_2}^2+\|\h_3\|_{\ell_2}^2=\|H\|_F^2$ we end up with:

\begin{lem}
\label{lem:psdweak1}
\begin{equation}
f(H,S_{wp})\leq \sqrt{ \| H\|_F^2-\sum_{i=1}^ch_{2,i}^2-\frac{(s(\h_1)+s(\h_2)-\sum_{i=1}^c h_{2,i})^2}{n-c}}
\end{equation}
for any $0\leq c\leq n-r$ such that $s(\h_1)+s(\h_2)-\sum_{i=1}^c h_{2,i}\geq (n-c)h_{2,c}$. If there is no such $c$ then $f(H,S_{wp})\leq \|H\|_F$
\end{lem}

Based on (\ref{lem:psdweak1}), for probabilistic analysis we'll use the following lemma:
\begin{lem}
\label{lem:psdweak2}
Let $H$ be chosen from $\HH(n)$ and $\h_1,\h_2,\h_3$ are vectors as described previously. Then we have: $f(H,S_{wp})\leq B_{wp}$ where
\begin{align}
&\nonumber B_{wp}= \| H\|_{F}~~~\text{if }g(H,c_{wp})\leq 0\\
&\label{eq:psdw1}B_{wp}=\sqrt{ \| H\|_F^2-\sum_{i=1}^{c_{wp}}h_{2,i}^2-\frac{(s(\h_1)+s(\h_2)-\sum_{i=1}^{c_{wp}} h_{2,i})^2}{n-{c_{wp}}}}~~~\text{else}
\end{align}
where $g(H,c)=\frac{s(\h_1)+s(\h_2)-\sum_{i=1}^c h_{2,i}}{n-c}-h_{2,c}$ and $c_{wp}=\delta_{wp} n(1-\beta)$ is a $c\leq n(1-\beta)$ such that
\begin{align}
\label{eq:csweakp}
c_{wp}&\text{ is solution of }~~~(1-\eps)\frac{\E [s(\h_1)+s(\h_2)-\sum_{i=1}^c h_{2,i}]}{\sqrt{n(1-\beta)}(n-c)}=F_s^{-1}\left(\frac{(1+\eps)c}{n(1-\beta)}\right)
\end{align}
where $\eps>0$ can be arbitrarily small.
\end{lem}
Note that $c_{wp}>0$ for any $\beta<1$, since we have $\E[s(\h_1)]=0$ and $\E[s(\h_2)]>0$.

\subsubsection{Probabilistic Analysis for $\E[B_{wp}]$}
Similar to probabilistic analysis for previous cases, we can use Lemma (\ref{lem:eiglem}) to show Lipschitzness of the function $s(\h_1)+s(\h_2)-\sum_{i=1}^c h_{2,i}$. Proof follows the exact same steps of Lemma (\ref{lem:seclip}). Then, using this and Lemmas \ref{lem:eiglip} and \ref{lem:gauslip}; we can conclude that $\Prob(g(H,c_{wp})\leq 0)$ decays to $0$ exponentially fast ($\exp(-O(n))$). As a result for $\E[B_{wp}]$ we have following upper bound by taking expectation of righthand side of (\ref{eq:psdw1}):
\beq
\E[B_{wp}]\leq n\sqrt{1-(1-\beta)^2(\gamma_{2,s}(1)-\gamma_{2,s}(1-\delta_{wp}))-\frac{(1-\beta)^3\gamma_s(1-\delta_{wp})^2}{1-(1-\beta)\delta_{wp}}+o(1)}+o(1)
\eeq
Then using Theorem (\ref{thm: ETM}) and $\gamma_{2,s}(1)=1$, we can conclude that:
\begin{thm}
\label{thm:psdwmain}
\beq
\mu>1-(1-\beta)^2(1-\gamma_{2,s}(1-\delta_{wp}))-\frac{(1-\beta)^3\gamma_s(1-\delta_{wp})^2}{1-(1-\beta)\delta_{wp}}
\eeq
is sufficient sampling rate for $\beta$ to be PSD weak threshold of Gaussian operator $\A:\R^{n\times n}\rightarrow\R^{\mu n(n+1)/2}$. Here, due to (\ref{eq:csweakp}), $\delta_{wp}$ is solution of:
\beq
(1-\eps)\frac{(1-\beta)^{3/2}\gamma_s(1-\delta)}{1-(1-\beta)\delta}=\sqrt{1-\beta}F_s^{-1}((1+\eps)\delta)
\eeq
Corresponding number of samples will be $m=\mu n(n+1)/2$. Also $F_s(\cdot)$ is the c.d.f. of the semicircle distribution defined previously.
\end{thm}
Remember that we choose smallest such $\mu$ to plot the curves. Result is given in Figure (\ref{fig:psd_thresholds}) as "Trace Minimization Weak".

\subsubsection{Alternative Analysis}
One can also directly analyze program (\ref{eq:optim4}) because it is exactly same as (\ref{eq:optim3}). It would give
\begin{lem}
\label{lem:psdweak3}
\begin{equation}
\label{eq:alter}
f(H,S_{wp})\leq \sqrt{ \|\h_1\|_{\ell_2}^2+\|\h_3\|_{\ell_2}^2+\|\h_4\|_{\ell_2}^2 -\sum_{i=1}^ch_{4,i}^2-\frac{(s(\h_1)+s(\h_4)-\sum_{i=1}^c h_{4,i})^2}{t-c}}
\end{equation}
for any $0\leq c\leq t-r$ such that $s(\h_1)+s(\h_4)-\sum_{i=1}^c h_{4,i}\geq (t-c)h_{4,c}$ where $t=t(\h_1,\h_4)=r+\eta_+(H_{2,2})$ is sum of dimensions of vectors $\h_1$ and $\h_4$. If there is no such $c$ then $f(H,S_{wp})\leq \|H\|_F$
\end{lem}
Note that $t$ is not deterministic however when $H$ is drawn from $\HH(n)$ we have $\E[t]=n(\beta+\frac{1-\beta}{2})=n\frac{1+\beta}{2}$ because clearly half of the eigenvalues of $H_{2,2}$ is positive in expectation. Furthermore from Lemmas \ref{lem:gauslip}, \ref{lem:eiglip}, it immediately follows that $t$ will concentrate around $\E[t]$ because for any $\eps$ we can write
\beq
\Prob(t-\frac{n(1+\beta)}{2}>n(1-\beta)\eps)= \Prob(\la_{n(1-\beta)(1/2+\eps)}(H_{2,2})\geq 0)\leq \exp(-\frac{n(1-\beta)}{2}F_s^{-1}(1/2-\eps)^2)
\eeq
(note that $F_s^{-1}(1/2)=0$) Then asymptotically $t/n$ will be approximately constant. Consequently we can probabilistically analyze Lemma (\ref{lem:psdweak3}) in a similar manner to previous cases however we have to deal with more details. At the end, we get the following which comes from asymptotic expectation of righthand side of (\ref{eq:alter}):
\begin{lem}
\beq
\mu>1-(1-\beta)^2(1-\frac{\gamma_2(1-\delta_{wp})}{2})-\frac{(1-\beta)^3\gamma(1-\delta_{wp})^2}{2(1+\beta-(1-\beta)\delta_{wp})}
\eeq
is sufficient sampling rate for $\beta$ to be PSD weak threshold of Gaussian operator $\A:\R^{n\times n}\rightarrow\R^{\mu n(n+1)/2}$. Here $\delta_{wp}$ is solution of
\beq
\frac{(1-\beta)^{3/2}\gamma(1-\delta)}{1+\beta-(1-\beta)\delta}=\sqrt{1-\beta}F^{-1}(\delta)
\eeq
\end{lem}
This formulation is nicer, since we did not use additional functions $F_s,\gamma_s,\gamma_{2,s}$.

\subsection{PSD Strong Threshold}
Now we'll analyze strong threshold for positive semidefinite matrices. 
\begin{PSDStr}
We say $\beta$ is a PSD strong threshold for Gaussian operator $\A:\R^{n\times n}\rightarrow\R^m$ ($m={\mu n(n+1)/2}$), if $A$ satisfies the following condition, asymptotically with probability $1$:

Any positive semidefinite matrix $X$ of rank at most $\beta n$ can be recovered from measurements $\A(X)$ via (\ref{eq: trace min}).
\end{PSDStr}

\begin{lem}
Any $X\in\PS_+^n$ of rank at most $r$ can be recovered from measurements $\A(X)$ via (\ref{eq: trace min}) if and only if any $W\in\N(\A)$ satisfies one of the following properties:
\begin{align}
&W~\text{is not hermitian}~~~\text{or}\\
&\text{trace}(W)>0~~~~~~~~~~~\text{or}\\
&\eta_-(W)>r
\end{align}
\begin{proof}
If one of the first two holds, then either $X+W$ is not PSD or $\text{trace}(X+W)>\text{trace}(X)$ so $X+W$ can not be a minimizer. On the other hand if third property holds then from Lemma (\ref{lem:inertia}) we find:
\beq
\eta_-(X+W)\geq \eta_-(W)-\eta_+(X)\geq r+1-\text{rank}(X)>0
\eeq
hence $X+W$ is not PSD. So $X$ will be unique minimizer if this is true for all $W$.

Conversely if there is a $W$ which satisfies none of the properties, then write $W=W_+-W_-$ where $W_+,W_-$ is PSD. Let $X=W_-$. Clearly $\text{rank}(X)=\eta_-(W)\leq r$ however $X+W$ is PSD and $\text{trace}(X+W)\leq \text{trace}(X)$ hence $X$ is not unique minimizer.
\end{proof}
\end{lem}

Now we'll analyze this condition similar to weak threshold for PSD matrices. $r=\beta n$, $\mu=n(n+1)/2$. Let $S_{sp}$ be the set of Hermitian matrices $W$ so that $\text{trace}(W)\leq0$, $\eta_-(W)\leq r$ and $\|W\|_F=1$. We don't want $S_{sp}$ to intersect with $\N_s(\A)$. Similar to previous analysis we need to calculate $\E[\sup_{W\in S_{sp}}\left<H,W\right>]$ to find the minimum sampling rate which ensures that intersection will be empty with high proabability.

Suppose $H\in\PS^n$ is fixed. Then let us calculate an upper bound $B_{sp}$ of $f(H,S_{sp})=\sup_{W\in S_{sp}}\left<H,W\right>$. If $\eta_-(H)<r$, we'll set $B_{sp}=\|H\|_F$ which is the obvious bound.

Otherwise from Lemma (\ref{lem:psdsimple}) we find:
\beq
\left<H,W\right>\leq \left<H_-,W_-\right>+\left<H_+,W_+\right>
\eeq
Let $c_+=\min\{\eta_+(H),\eta_+(W)\}$ and $c_-=\min\{\eta_-(H),\eta_-(W)\}$ Then from Lemma (\ref{lem:neumann}):
\beq
\left<H,W\right>\leq u(H,W):=\sum_{i=1}^{c_+}\la_i(H_+)\la_i(W_+)+\sum_{i=1}^{c_-}\la_i(H_-)\la_i(W_-)
\eeq
To upper bound $f(H,S_{sp})$ let us maximize $u(H,W)$ over $S_{sp}$. Let $W\in S_{sp}$ then we have $\eta_-(H)\geq r\geq \eta_-(W)$.


%
Let $\h_1,\w_1\in\R^{\eta_+(H)}$ be vectors increasingly ordered largest $\eta_+(H)$ eigenvalues of $H_+,W_+$ respectively. Similarly $\h_2,\w_2\in\R^{r}$ be vectors of increasingly ordered largest $r$ eigenvalues of $H_-,W_-$. Since $\eta_+(H)\geq c_+$ and $r\geq c_-$ we can write:
\beq
\h_1^T\w_1+\h_2^T\w_2=u(H,W)
\eeq
Note that $\w_1,\w_2$ has to satisfy: $\w_1,\w_2\succeq 0$, $s(\w_1)\leq s(\w_2)$ and $\|\w_1\|_{\ell_2}^2+\|\w_2\|_{\ell_2}^2\leq\|W\|_F^2\leq 1$. Middle one is due to $s(\w_1)\leq \text{trace}(W_+)$ and $s(\w_2)=\text{trace}(W_-)$ and $\text{trace}(W)= \text{trace}(W_+)- \text{trace}(W_-)\leq 0$. Also we have $\h_1,\h_2\succeq 0$. Then following optimization program will give $\sup_{W_\in S_{sp}} u(H,W)$
\begin{align}
&\max_{\y_1,\y_2}~\h_1^T\y_1+\h_2^T\y_2\\
&\nonumber \text{subject to}\\
&\nonumber~~~~\y_1,\y_2\succeq 0\\
&\nonumber~~~~s(\y_1)\leq s(\y_2)\\
&\nonumber~~~~\|\y_1\|_{\ell_2}^2+\|\y_2\|_{\ell_2}^2\leq 1
\end{align}
Note that this is exactly same as program (\ref{eq:optim1}). As a result we can write the following Lemma:
\begin{lem}
\label{lem:psdstr1}
Let $t=\eta_+(H)+r$.
\begin{align}
&f(H,S_{sp})\leq \|H\|_F~~\text{if}~~\eta_-(H)<r~\text{or}~s(\h_1)\leq s(\h_2)\\
\label{eq:only}
&f(H,S_{sp})\leq \sqrt{\|\h_1\|_\lt^2+\|\h_2\|_\lt^2-\sum_{i=1}^c h_{1,i}^2-\frac{(s(\h_1)-s(\h_2)-\sum_{i=1}^ch_{1,i})^2}{t-c}}~~\text{else}
\end{align}
where $c\leq \eta_+(H)$ such that $s(\h_1)-s(\h_2)-\sum_{i=1}^ch_{1,i}\geq (t-c)h_{1,c}$
\end{lem}
We'll not give the detailed ETM analysis for this case, as it requires more meticulous analysis. However for any $r=\beta n$ with $\beta<1/2$ it is easy to show that, when $H$ is chosen from $\HH(n)$, we'll have
\beq
\Prob(\eta_-(H)<\beta n~\text{or}~s(\h_1)\leq s(\h_2))\rightarrow 0
\eeq
exponentially fast with $n$. The reason is that $\E[\eta_-(H)]=n/2>\beta n$ and $\E[s(\h_1)- s(\h_2)]=\frac{n^{3/2}}{2}(\gamma(1)-\gamma(2\beta))>0$. Similar to previous cases, using Lipschitzness of the functions and Gaussianity of $H$ will yield the result. Then essentially we need to analyze (\ref{eq:only}). Using exactly same arguments we can also show, $\E[t]=n(\beta+1/2)$ and $t$ will concentrate around its mean (as $n\rightarrow\infty$). As a result, except minor details, probabilistic analysis becomes similar to the ones before (i.e. where $t$ is constant). At the end we get:
\begin{thm}
\label{eq:psdstr2}
If $\beta \geq 1/2$ then $\mu=1$. When $\beta<1/2$ we have:
\beq
\mu>\frac{1}{2}\left[\gamma_2(1-\delta_{sp})+\gamma_2(2\beta)-\frac{(\gamma(1-\delta_{sp})-\gamma(2\beta))^2}{2\beta+1-\delta_{sp}}\right]
\eeq
is a sufficient sampling rate for $\beta$ to be PSD strong threshold of Gaussian operator $\A:\R^{n\times n}\rightarrow\R^{\mu n(n+1)/2}$. Here $\delta_{sp}$ is solution of
\beq
\frac{\gamma(1-\delta)-\gamma(2\beta)}{2\beta+1-\delta}=F^{-1}(\delta)
\eeq
\end{thm}
Result is given in Figure \ref{fig:psd_thresholds} as "Trace Minimization Strong".
 
\subsection{Uniqueness Results}
In this part, we'll state the conditions and results for the unique PSD solution to the measurements without proof. They follow immediately from slight modifications of previous analysis.
\subsubsection{Weak Uniqueness}
\begin{UniWeak}
Let $\A:\R^{n\times n}\rightarrow \R^m$ be a random Gaussian operator and let $X$ be an arbitrary PSD matrix with $\text{rank}(X)=\beta n$. We say that $\beta$ is a uniqueness weak threshold if with high probability this particular matrix $X$ can be recovered from measurements $\A(X)$ via program (\ref{eq: unique PSD}).
\end{UniWeak}
\begin{lem}
Let $X$ be a PSD matrix with $\text{rank}(X)=r$ and eigenvalue decomposition $U\Lambda U^T$ with $\Lambda\in\R^{r\times t}$. Then $X$ can be recovered via (\ref{eq: unique PSD}) if for all $W\in\N_s(\A)$, $\bar{U}^TW\bar{U}$ has a negative eigenvalue.
\end{lem}
\begin{lem}
Let $\A:\R^{n\times n}\rightarrow \R^{\mu n(n+1)/2}$ be a Gaussian operator. Then $\beta$ is a weak uniqueness threshold if
\beq
\mu>1-\frac{(1-\beta)^2}{2}
\eeq
\end{lem}

\subsubsection{Strong Uniqueness}

\begin{UniStr}
Let $\A:\R^{n\times n}\rightarrow \R^m$ be a random Gaussian operator. We say that $\beta$ is a uniqueness strong threshold if with high probability all PSD matrices $X$ with rank at most $\beta n$ can be recovered from their measurements $\A(X)$ via program (\ref{eq: unique PSD}).
\end{UniStr}
\begin{lem}
All PSD matrices of rank at most $r$ can be recovered via program (\ref{eq: unique PSD}) if and only if all $W\in\N_s(\A)$, has at least $r+1$ negative eigenvalue.
\end{lem}
\begin{lem}
Let $\A:\R^{n\times n}\rightarrow \R^{\mu n(n+1)/2}$ be a Gaussian operator. Then $\beta$ is a strong uniqueness threshold if
\begin{align}
&\mu=1~~~~~~~~~~~~~~~\text{if}~\beta\geq 0.5\\
&\mu=\frac{1+\gamma_2(2\beta)}{2}~~~\text{else}
\end{align}
\end{lem}

\begin{figure}[t]
\centering
  \includegraphics[width= 1\textwidth]{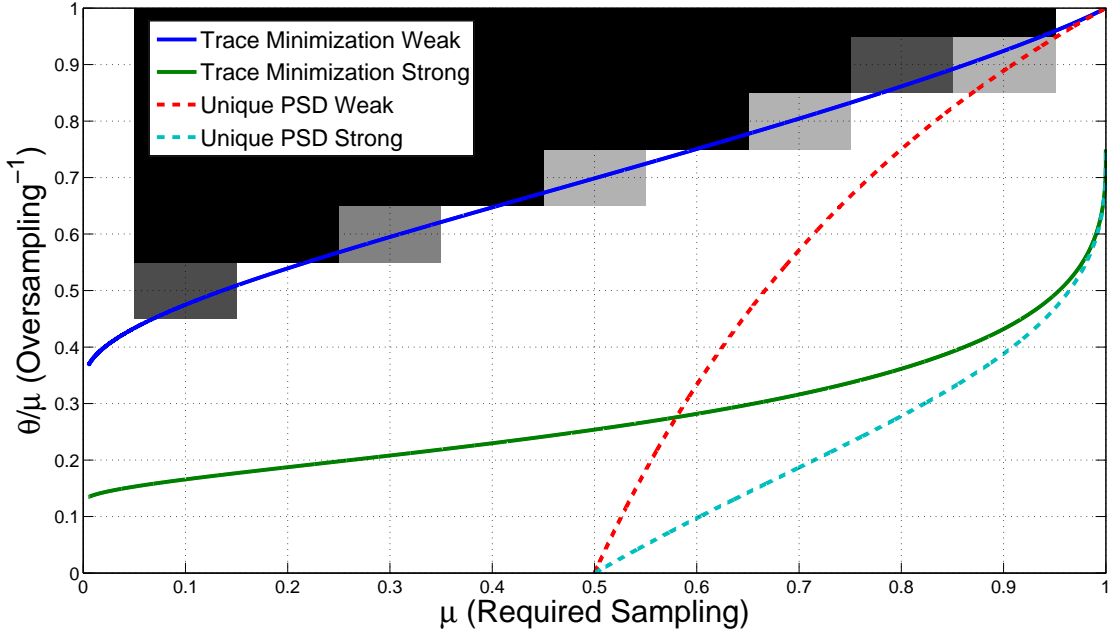}
 \caption{\scriptsize{Results for PSD matrices. Again $\text{oversampling}^{-1}$ ($\theta/\mu$) vs $\mu$ is plotted.} Simulations are done for $40\times 40$ matrices and program (\ref{eq: trace min}) is solved with Gaussian measurements. Although resolution is low (simulations are not fine), it is not hard to see that trace minimization weak threshold looks consistent with simulations. Black and white regions mean failure and success respectively.}
  \label{fig:psd_thresholds}
\end{figure}

Curves for weak and strong uniqueness thresholds are given in Figure (\ref{fig:psd_thresholds}) as "Unique PSD Weak/Strong".

\section{Discussion and Future Works}
In this work we classified the various types of matrix recovery, gave tight conditions for them and analyzed the conditions for Gaussian measurements to get better thresholds than the existing results of \cite{Recht_Xu_Hassibi} and \cite{arxiv}. It turns out that the thresholds of \cite{Recht_Xu_Hassibi,arxiv} actually corresponds to a special, suboptimal case of our analysis. In Lemmas \ref{lem:strpes1}, \ref{lem:secpes1}, \ref{lem:weak2} instead of choosing $\delta_s,\delta_{sec},\delta_w$ carefully if we just set them to $0$, we'll end up with results of \cite{Recht_Xu_Hassibi,arxiv}. This suggests that, although analysis of this paper is more tedious, it is strictly better than previous ones and also generalizes them.

Although we didn't do much argument about tightness of our results, actually most of the estimations and inequalities that are used for upper boundings are tight or asymptotically tight. In particular we believe our weak thresholds are exact, similar to the significant results of \cite{mihailo}. A key to the results of the paper is the fact that we have written down the null space conditions in their most transparent form. Essentially, the null space vectors of compressed sensing are replaced by the singular values of the null space matrix in NNM. This allowed us to use the approach of \cite{mihailo} directly. This furthermore suggests that the NNM problem is a generalization of compressed sensing and the two problems are very similar in nature.

Our simulation results support our belief that our weak thresholds are tight. Also simulation results and theoretical curves suggest that at most 3 times of oversampling is necessary for weak recovery for any $0\leq \beta\leq 1$, and around 8 times is required for strong. This is important as it means one can solve the RM problem via convex optimization with a very small sampling cost. Furthermore, although our results are in the asymptotic case ($r=\beta n$), theory and simulation fits almost perfectly even for a relatively small matrix of size $40\times 40$. This suggests that actually concentration of measure happens pretty quickly.

It would be interesting to calculate the limiting case of $\beta/\mu$ as $\beta\rightarrow 0$, to get an estimate of the minimum required oversampling when the rank is small. Secondly, we believe it might be possible to employ these methods not only for the linear region where rank $r=\beta n$ but for any case such as $r=O(1)$ or $r=O(log(n))$. Such a study might give a small ($\approx 3,4$) minimum oversampling rate. Although recent results of \cite{candes_last} showed we need only $O(rn)$ samples for recovery (which is minimal), the constant is not known.

Finally, our result suggests a significant performance difference between trace minimization and unique solution in the special case of PSD matrices. Although we'll not argue the reason here, it is actually quite intuitive. Uniqueness results suggests that one needs to sample at least half of the entries ($n(n+1)/4$ samples for PSD) to make sure that positive semidefinite solution is unique. Clearly such a result might be interesting to know but it is not useful at all.

Our first aim will be verifying our tightness claims. In order to do this, one needs to investigate the work of \cite{Gordan2} better and to come up with the conditions on the "mesh" where (\ref{thm: ETM}) is tight.

{\bf{Comment: }}Overall, the results of \cite{mihailo} established a powerful way to analyze some important questions in low rank matrix recovery.

\end{document}